\documentclass[11pt]{amsart}

\usepackage{amsmath,amssymb,amsthm}
\usepackage[margin=1.12in]{geometry}
\usepackage{hyperref}
\hypersetup{hidelinks}

\theoremstyle{plain}
\newtheorem{theorem}{Theorem}[section]
\newtheorem{proposition}[theorem]{Proposition}
\newtheorem{lemma}[theorem]{Lemma}
\newtheorem{corollary}[theorem]{Corollary}
\theoremstyle{definition}

\newtheorem{conjecture}[theorem]{Conjecture}

\newcommand{\e}{\mathrm e}
\newcommand{\dd}{\,\mathrm d}
\newcommand{\abs}[1]{\lvert #1\rvert}

\begin{document}

\title{Monotonicity of Erd\H{o}s sums}
\author{Jaideep Sai Padhi}
\address{Department of Computer Science, Purdue University,
West Lafayette, IN 47907, USA}
\email{jpadhi@purdue.edu}
\date{}

\begin{abstract}
For $k\geq1$ let $f_k=\sum_{\Omega(n)=k}1/(n\log n)$, where $\Omega$ counts
prime factors with multiplicity, and let $g_k$ be the same sum restricted to
odd $n$. Lichtman showed that $f_k\to1$, with a global minimum at $k=6$, and
Gorodetsky, Lichtman and Wong that
$f_k=1-2^{-k}\bigl(a_1k^2+O(k\log(k+1))\bigr)$. We prove that $f_k<f_{k+1}$
for every $k\geq6$, and that $g_k>g_{k+1}$ for every $k\geq1$. The two
statements come from one analytic distinction. Writing
$\sum_kP_k(s)z^k=\zeta(s)^zH(s,z)$, at $s=1$ the point $z=2$ lies on the first
singular circle of $H$, and strictly inside the first singular circle of the
product with the prime $2$ removed. We show that this distinction decides the
direction in general. If the complement of a set of primes $\mathcal P$ is
thin and its least element is at least $3$, the corresponding sums eventually
decrease, so the second theorem is one case of a general phenomenon. The
inequalities used at the last step of the first
theorem are certified by interval arithmetic.
\end{abstract}

\maketitle

\section{Introduction}

For $k\geq1$ put
\[
  f_k=\sum_{\Omega(n)=k}\frac1{n\log n}.
\]
Erd\H{o}s \cite{Erdos} proved that the corresponding sum is bounded uniformly
over primitive sets. In particular each $f_k$ converges. Zhang \cite{Zhang}
proved $f_k<f_1$ for $k>1$, and Bayless, Kinlaw and Klyve \cite{BKK} proved
$f_2>f_3$. Banks and Martin \cite{BanksMartin} conjectured that
\begin{equation}\label{eq:BM}
  f_1>f_2>f_3>\cdots.
\end{equation}

Lichtman \cite{Lichtman19} disproved \eqref{eq:BM}, showing that $f_k\to1$ and
that the sequence attains a global minimum at $k=6$. Gorodetsky, Lichtman and
Wong \cite{GLW} then established
\begin{equation}\label{eq:GLW-asymptotic}
  f_k=1-2^{-k}\bigl(a_1k^2+O(k\log(k+1))\bigr),
  \qquad a_1=\frac{d\log2}{4}>0,
\end{equation}
where
\begin{equation}\label{eq:d-def}
  d=\frac14\prod_{p\geq3}\left(1-\frac2p\right)^{-1}
                    \left(1-\frac1p\right)^2.
\end{equation}
Neither of these yields monotonicity beyond the minimum. A global minimum at
$k=6$ leaves open whether $f_k<f_{k+1}$ holds at every index above it, and the
implicit constant in \eqref{eq:GLW-asymptotic} does not make explicit a range
on which the increase can be read off. Our first theorem gives that
monotonicity in full.

\begin{theorem}\label{thm:main-full}
For every $k\geq6$ we have $f_k<f_{k+1}$.
\end{theorem}

Let
\[
  g_k
   =\sum_{\substack{n\ \mathrm{odd}\\ \Omega(n)=k}}\frac1{n\log n}.
\]
The second theorem shows that deleting the prime $2$ reverses the direction of
monotonicity at every index.

\begin{theorem}\label{thm:main-odd}
For every $k\geq1$ we have $g_k>g_{k+1}$, and $g_k\to\tfrac12$.
\end{theorem}

The proof of Theorem~\ref{thm:main-odd} gives the more precise expansion
\begin{equation}\label{eq:odd-intro-expansion}
  g_k=\frac12+\beta2^{-k}+R_k,
  \qquad \abs{R_k}\leq4900(k+3)^3 3^{-k}\ \ (k\geq59),
\end{equation}
where $\gamma=0.5772156649\ldots$ is Euler's constant and
\begin{equation}\label{eq:beta-intro}
  \beta=d\left(\gamma+\log2-
  \sum_{p\geq3}\frac{\log p}{(p-1)(p-2)}\right)>0.
\end{equation}
The $2^{-k}$ term comes from evaluating a regularized Euler product at $z=2$.
For odd integers the first pole is at $z=3$, so this evaluation takes place
strictly inside the disc of analyticity. For the unrestricted sums the first
pole is at $z=2$ itself. The boundary contribution then acquires the factor
$k^2$ in \eqref{eq:GLW-asymptotic}.

This contrast is not an accident of the two cases. In
Section~\ref{sec:general-principle} we make it a general statement. Let
$\mathcal P$ be a set of primes with least element $q\geq3$, and suppose its
complement $\mathcal Q$ is thin, in the sense that
$\sum_{p\in\mathcal Q}(\log p)/p<\infty$. Then
Proposition~\ref{prop:general-interior} gives
\[
  f_k^{\mathcal P}
   =\prod_{p\in\mathcal Q}\Bigl(1-\frac1p\Bigr)+c_{\mathcal P}2^{-k}+o(2^{-k}),
  \qquad
  c_{\mathcal P}=H_{\mathcal P}(1,2)
   \bigl(\gamma+\tfrac12\partial_s\log H_{\mathcal P}(1,2)\bigr),
\]
and Corollary~\ref{cor:general-positivity} gives $c_{\mathcal P}>0$, so every
such family eventually decreases. Thus Theorem~\ref{thm:main-odd} is the case
$\mathcal Q=\{2\}$ of a phenomenon that depends only on the position of $z=2$
relative to the first singular circle. Both hypotheses are used. Thinness does
not follow from $q\geq3$: for $\mathcal P=\{p\equiv1\bmod4\}$ the leading
product is $0$, and the proposition says nothing.

The remaining case is $2\in\mathcal P$, where the evaluation sits on the
singular circle itself. That case is open, and
Conjecture~\ref{conj:general} records what it would give; we comment on it at
the end of Section~\ref{sec:general-principle} but do not treat it here.

The proof has two kinds of input. Identities, coefficient reductions and
inequalities are proved in the text. Bounds produced by interval arithmetic
are isolated as certified propositions, so that a phrase of the form ``the
computation gives'' never stands in for an analytic step: the exact inequality
the computation supplies is stated before it is used. What the interval
arithmetic certifies is the computation itself. Transcribing an enclosure into
a numeral printed here, and combining several such numerals, is ordinary
writing, and we check it by script rather than by ball arithmetic.

\section{Generating functions and the analytic boundary}

For $s>1$ define
\[
  P_k(s)=\sum_{\Omega(n)=k}n^{-s},
  \qquad
  P_k^{\mathrm{odd}}(s)=
  \sum_{\substack{n\ \mathrm{odd}\\\Omega(n)=k}}n^{-s}.
\]

\begin{lemma}\label{lem:integral-representation}
For every $k\geq1$,
\begin{equation}\label{eq:integral-representation}
  f_k=\int_1^\infty P_k(s)\dd s,
  \qquad
  g_k=\int_1^\infty P_k^{\mathrm{odd}}(s)\dd s.
\end{equation}
\end{lemma}

\begin{proof}
The first identity is \cite[(3.1)]{Lichtman19}: one has
$\int_1^\infty n^{-s}\dd s=1/(n\log n)$ for $n\geq2$, and the summands are
non-negative, so Tonelli's theorem permits the interchange. The same
calculation with the sum restricted to odd $n$ gives the second.
\end{proof}

The generating functions are, by \cite[Prop.~3.1]{Lichtman19},
\begin{equation}\label{eq:euler-products}
  \sum_{k\geq0}P_k(s)z^k=\prod_p(1-zp^{-s})^{-1},
  \qquad
  \sum_{k\geq0}P_k^{\mathrm{odd}}(s)z^k
      =\prod_{p\geq3}(1-zp^{-s})^{-1}.
\end{equation}
Write these products as $\zeta(s)^zH(s,z)$ and
$\zeta(s)^zH_{\mathrm o}(s,z)$, where
\begin{align}
  H(s,z)&=\prod_p(1-zp^{-s})^{-1}(1-p^{-s})^z,
       \label{eq:H-full}\\
  H_{\mathrm o}(s,z)&=(1-2^{-s})^z
       \prod_{p\geq3}(1-zp^{-s})^{-1}(1-p^{-s})^z.
       \label{eq:H-odd}
\end{align}

\begin{proposition}\label{prop:analyticity}
Let $\sigma=\operatorname{Re}s$. Then $H$ is analytic in the two complex
variables $(s,z)$ on $\{\sigma>1/2,\ \abs z<2^{\sigma}\}$, and
$H_{\mathrm o}$ on $\{\sigma>1/2,\ \abs z<3^{\sigma}\}$. In particular, for
fixed real $s>1/2$ the function $H(s,\cdot)$ is analytic in $\abs z<2^s$ and
$H_{\mathrm o}(s,\cdot)$ in $\abs z<3^s$; each
$\partial_s^rH_{\mathrm o}(s,\cdot)$ is again analytic in
$\abs z<3^s$, and
\begin{equation}\label{eq:coeff-commute}
  \frac{\dd^r}{\dd s^r}[z^j]H_{\mathrm o}(s,z)
   =[z^j]\partial_s^rH_{\mathrm o}(s,z);
\end{equation}
the same holds for $H$ in $\abs z<2^s$.
\end{proposition}

\begin{proof}
It is essential that $s$ be taken complex here: uniform convergence of
real-variable functions says nothing about the derivative of the limit, so the
real-variable argument would give analyticity in $z$ alone.

Let $x=p^{-s}$, so that $\abs x=p^{-\sigma}<2^{-1/2}$ on the domain. The
factor $(1-x)^z$ needs a branch, and the principal one serves throughout,
since $1-x$ lies in the disc of radius $2^{-1/2}<1$ about $1$. Fix a compact
subset of the domain, so that $\sigma\geq s_0>1/2$ and $\abs z\leq R<2^{s_0}$
for \eqref{eq:H-full}, or $R<3^{s_0}$ with $p\geq3$ for \eqref{eq:H-odd}; then
$\abs{zx}<1$ for every $p$ occurring. The logarithm of the $p$th regularized
factor is
\begin{align}
 -\log(1-zx)+z\log(1-x)
 &=\sum_{m\geq1}\frac{z^m-z}{m}x^m
  =\sum_{m\geq2}\frac{z^m-z}{m}x^m
  =O(p^{-2\sigma}),                       \label{eq:reg-factor-log}
\end{align}
the $m=1$ terms cancelling. The estimate is uniform on the compact set and
depends on $s$ only through $\sigma$. Each $p^{-s}=\e^{-s\log p}$ is entire,
so every regularized factor is analytic in $(s,z)$ jointly wherever its
denominator does not vanish, which on each of the two domains it does not: for
\eqref{eq:H-full} one has $\abs z<2^{\sigma}\leq p^{\sigma}=\abs{p^{s}}$ for
every prime $p$, and for \eqref{eq:H-odd} $\abs z<3^{\sigma}\leq\abs{p^{s}}$
for every $p\geq3$. As $\sum_pp^{-2s_0}<\infty$, the sum of the logarithms
converges absolutely and locally uniformly, so the product converges locally
uniformly on the stated domain in $\mathbb C^2$ to a non-vanishing limit, and
a locally uniform limit of analytic functions of several complex variables is
analytic. This is the first assertion. The factors have poles at $z=p^s$, of
modulus $p^{\sigma}$, so the nearest is at modulus $2^{\sigma}$ for
\eqref{eq:H-full} and $3^{\sigma}$ for \eqref{eq:H-odd}, which is why the
domains are as stated. Finally $[z^j]$ is a
Cauchy integral over a circle that may be held fixed on a compact subset, so
it commutes with $\dd/\dd s$, which is \eqref{eq:coeff-commute}.
\end{proof}

At $s=1$ the point $z=2$ lies on the first singular circle for $H$ and lies
strictly inside the first singular circle for $H_{\mathrm o}$. This is the
analytic distinction used below.

We shall also need the local change of variables used in the odd case.

\begin{lemma}\label{lem:zeta-coordinate}
Let $u=\log\zeta(s)$ for $s>1$, and let $s=s(u)$ be the inverse function. Put
$\varepsilon(u)=s(u)-1$ and $J(u)=-s'(u)$. As $u\to\infty$,
\begin{equation}\label{eq:zeta-coordinate}
  \varepsilon(u)=\e^{-u}+\gamma\e^{-2u}+O(\e^{-3u}),
  \qquad
  J(u)=\e^{-u}+2\gamma\e^{-2u}+O(\e^{-3u}).
\end{equation}
\end{lemma}

\begin{proof}
The point of the argument is that the expansions below are convergent power
series and not merely asymptotic, which is what allows the second of them to be
differentiated term by term.

By the Laurent expansion of $\zeta$ at its simple pole
\cite[Ch.~1]{MV}, the function $\varepsilon\mapsto\varepsilon\,
\zeta(1+\varepsilon)$ is analytic near $\varepsilon=0$ with expansion
$1+\gamma\varepsilon+O(\varepsilon^2)$. Hence
\[
 x=\frac1{\zeta(1+\varepsilon)}
  =\frac{\varepsilon}{\varepsilon\,\zeta(1+\varepsilon)}
\]
is analytic near $\varepsilon=0$, with $x(0)=0$ and
$x=\varepsilon-\gamma\varepsilon^2+O(\varepsilon^3)$, so in particular
$\dd x/\dd\varepsilon=1\neq0$ at the origin. By the inverse function theorem
for analytic functions, $\varepsilon$ is therefore an analytic function of $x$
near $x=0$, given by a convergent power series
\[
 \varepsilon=\psi(x)=x+\gamma x^2+O(x^3),
\]
which is the first formula in \eqref{eq:zeta-coordinate} once
$x=\e^{-u}$. Because $\psi$ is analytic rather than merely asymptotic, it may
be differentiated termwise: $\psi'(x)=1+2\gamma x+O(x^2)$. With
$x=\e^{-u}$ one has $\dd x/\dd u=-x$, so the chain rule gives
\[
  \varepsilon'(u)=\psi'(x)\frac{\dd x}{\dd u}=-x\bigl(1+2\gamma x+O(x^2)\bigr)
  =-x-2\gamma x^2+O(x^3),
\]
and $J(u)=-\varepsilon'(u)=x+2\gamma x^2+O(x^3)$.
\end{proof}

\section{The odd sums}\label{sec:odd}

Throughout this section write
\[
  h_j(s)=[z^j]H_{\mathrm o}(s,z),
  \qquad
  b_j=h_j'(1)+2\gamma h_j(1).
\]
Since
$\zeta(s)^z=\sum_{m\geq0}(\log\zeta(s))^m z^m/m!$, coefficient extraction in
\eqref{eq:euler-products} gives the exact identity
\begin{equation}\label{eq:odd-coefficient-extraction}
  P_k^{\mathrm{odd}}(s)
  =\sum_{j=0}^k h_j(s)\frac{(\log\zeta(s))^{k-j}}{(k-j)!}.
\end{equation}
Using $u=\log\zeta(s)$ in Lemma~\ref{lem:integral-representation}, and noting
that $s$ decreases from $\infty$ to $1$ as $u$ increases from $0$ to
$\infty$, we obtain
\begin{equation}\label{eq:odd-u-integral}
  g_k
  =\sum_{j=0}^k\int_0^\infty
       \frac{u^{k-j}}{(k-j)!}\,h_j(s(u))J(u)\dd u.
\end{equation}

\subsection{The first two terms}

\begin{proposition}\label{prop:odd-two-term}
With $h_j$ and $b_j$ as above and $E_j$ as in the proof below, put
\begin{equation}\label{eq:Rk-def}
 R_k=\sum_{j=0}^k\int_0^\infty\frac{u^{k-j}}{(k-j)!}E_j(u)\dd u
     -\sum_{j>k}h_j(1)
     -2^{-k-1}\sum_{j>k}b_j2^j.
\end{equation}
Then
\begin{equation}\label{eq:odd-two-term}
 g_k
 =H_{\mathrm o}(1,1)
  +2^{-k-1}\bigl(\partial_sH_{\mathrm o}(1,2)
                  +2\gamma H_{\mathrm o}(1,2)\bigr)+R_k.
\end{equation}
\end{proposition}

\begin{proof}
The idea is to write the integrand of \eqref{eq:odd-u-integral} as its first
two terms in $\e^{-u}$ plus an explicit error, and then to integrate each
piece exactly. Two expansions are at work and they must be kept apart: the
measure $J(u)\dd u$, which Lemma~\ref{lem:zeta-coordinate} expands, and the
coefficient $h_j(s(u))$, which Taylor's theorem expands about $s=1$. Their
product supplies the three error terms of \eqref{eq:odd-local-product}.

Let $x=\e^{-u}$. By Lemma~\ref{lem:zeta-coordinate},
$J(u)=x+2\gamma x^2+O(x^3)$ and
$\varepsilon(u)J(u)=x^2+O(x^3)$. Taylor's theorem at $s=1$ gives
\[
  h_j(1+\varepsilon)
   =h_j(1)+\varepsilon h_j'(1)
     +\frac{\varepsilon^2}{2}h_j''(1+\theta_{j,u}\varepsilon)
\]
for some $0<\theta_{j,u}<1$. Put
\[
 r_0(u)=J(u)-\e^{-u}-2\gamma\e^{-2u},
 \qquad
 r_1(u)=\varepsilon(u)J(u)-\e^{-2u}.
\]
Multiplication by $J(u)$ now gives the exact identity
\begin{equation}\label{eq:odd-local-product}
\begin{aligned}
 h_j(s(u))J(u)
  &=h_j(1)\e^{-u}+b_j\e^{-2u}+E_j(u),\\
 E_j(u)
  &=h_j(1)r_0(u)+h_j'(1)r_1(u)
    +\frac{\varepsilon(u)^2J(u)}2
       h_j''(1+\theta_{j,u}\varepsilon(u)).
\end{aligned}
\end{equation}
No error is suppressed in the notation $E_j$: its three summands are the
Jacobian error, the first-order coordinate error and the Taylor remainder.

Insert \eqref{eq:odd-local-product} into \eqref{eq:odd-u-integral}. Since
\[
  \int_0^\infty\frac{u^m}{m!}\e^{-ru}\dd u=r^{-m-1}
  \qquad(r>0),
\]
the first two terms become
\begin{align}
 \sum_{j=0}^k h_j(1)
   &=H_{\mathrm o}(1,1)-\sum_{j>k}h_j(1),
   \label{eq:first-coefficient-tail}\\
 \sum_{j=0}^k b_j2^{-(k-j+1)}
   &=2^{-k-1}\bigl(\partial_sH_{\mathrm o}(1,2)
                    +2\gamma H_{\mathrm o}(1,2)\bigr)
     -2^{-k-1}\sum_{j>k}b_j2^j.
   \label{eq:second-coefficient-tail}
\end{align}
Substituting \eqref{eq:first-coefficient-tail} and
\eqref{eq:second-coefficient-tail}, and collecting the $E_j$ integrals with the
two tails, is exactly \eqref{eq:Rk-def}, and this gives
\eqref{eq:odd-two-term}. Note that \eqref{eq:Rk-def} is an identity, so every
later estimate is an estimate of those three explicit quantities.
\end{proof}

\begin{lemma}\label{lem:odd-constants}
One has
\begin{equation}\label{eq:odd-H-values}
  H_{\mathrm o}(1,1)=\frac12,
  \qquad
  H_{\mathrm o}(1,2)=d,
\end{equation}
and
\begin{equation}\label{eq:odd-log-derivative}
  \partial_s\log H_{\mathrm o}(1,2)
   =2\log2-2\sum_{p\geq3}\frac{\log p}{(p-1)(p-2)}.
\end{equation}
\end{lemma}

\begin{proof}
The first two values are immediate from \eqref{eq:H-odd}: at $s=z=1$ every
factor in the product over $p\geq3$ cancels, leaving
$H_{\mathrm o}(1,1)=1-2^{-1}=\tfrac12$, and at $(s,z)=(1,2)$ the product is
\eqref{eq:d-def}, so $H_{\mathrm o}(1,2)=d$.

For the derivative, take logarithms:
\[
 \log H_{\mathrm o}(s,2)
 =2\log(1-2^{-s})
  +\sum_{p\geq3}
     \left[-\log(1-2p^{-s})+2\log(1-p^{-s})\right].
\]
The differentiated $p$th summand is $O((\log p)p^{-2s})$ locally uniformly
near $s=1$, so termwise differentiation is permitted. At $s=1$ the first term
contributes $2\log2$, while the $p$th summand contributes
\begin{align*}
 -\frac{2\log p}{p-2}+\frac{2\log p}{p-1}
 =-\frac{2\log p}{(p-1)(p-2)}.
\end{align*}
Summing over $p\geq3$ proves \eqref{eq:odd-log-derivative}.
\end{proof}

Combining Proposition~\ref{prop:odd-two-term} and
Lemma~\ref{lem:odd-constants} gives
\begin{equation}\label{eq:odd-expansion}
  g_k=\frac12+\beta2^{-k}+R_k,
  \qquad
  \beta=d\left(\gamma+\log2-
  \sum_{p\geq3}\frac{\log p}{(p-1)(p-2)}\right).
\end{equation}

\subsection{The remainder}

Three weighted $\ell^1$ norms of coefficients carry all of the interval
input for the odd sums. Throughout, $s(3)=s(u)\rvert_{u=3}=1.0512\ldots$ in the
coordinate of Lemma~\ref{lem:zeta-coordinate}, that is, the value of $s$ at
which $\log\zeta(s)=3$. Put
\begin{equation}\label{eq:sigma-def}
\begin{aligned}
 \Sigma_0&=\sum_{m\geq0}3^m
   \abs{[z^m](1-z/3)H_{\mathrm o}(1,z)},\\
 \Sigma_1&=\sum_{m\geq0}3^m
   \abs{[z^m](1-z/3)^2\partial_sH_{\mathrm o}(1,z)},\\
 \Sigma_2&=\sup_{1\leq s\leq s(3)}
   \sum_{m\geq0}3^m
   \abs{[z^m](1-z/3)^3\partial_s^2H_{\mathrm o}(s,z)}.
\end{aligned}
\end{equation}
Finiteness of $\Sigma_0$ and $\Sigma_1$ is analytic, and the weights $r=1,2$ on
the factors $(1-z/3)^r$ are what makes it so.
Proposition~\ref{prop:analyticity} gives $H_{\mathrm o}(1,\cdot)$ only on
$\abs z<3$, which by itself would not make a sum weighted by $3^m$ converge.
But at $s=1$ the nearest pole in $z$ is at $z=3$, contributed by the single
Euler factor $(1-z/3)^{-1}$, and after $r-1$ derivatives in $s$ that pole has
order $r$; multiplying by $(1-z/3)^r$ cancels it exactly, so both functions
continue analytically to $\abs z<5$, the next pole being at $z=5$. Their
coefficients decay like $5^{-m}$ and the weighted sums converge geometrically.

For $\Sigma_2$ the same device applies only at the left endpoint. At $s>1$ the
pole sits at $z=3^s$, which $(1-z/3)^3$ does not cancel, and the inner sum is
finite because the radius of convergence is then $3^s>3$, a margin that
shrinks to nothing as $s\to1^+$, so this does not by itself bound the
supremum. Each inner sum is finite for every $s$ in the range, by cancellation
at $s=1$ and by $3/3^s<1$ for $s>1$; the finiteness of the supremum over
$[1,s(3)]$ is part of what the interval computation below establishes, not a
separate analytic claim. For $\Sigma_0$ and $\Sigma_1$ only the numerical
values are certified.

Next, split the remainder of Proposition~\ref{prop:odd-two-term} at $u=3$:
with $E_j$ the explicit expression in \eqref{eq:odd-local-product}, put
\begin{equation}\label{eq:near-far-split}
 R_{k,\mathrm{far}}=\sum_{j=0}^k\int_3^\infty
   \frac{u^{k-j}}{(k-j)!}E_j(u)\dd u,
 \qquad
 R_{k,\mathrm{near}}=\sum_{j=0}^k\int_0^3
   \frac{u^{k-j}}{(k-j)!}E_j(u)\dd u,
\end{equation}
so that $R_k$ is $R_{k,\mathrm{far}}+R_{k,\mathrm{near}}$ together with the two
coefficient tails of
\eqref{eq:first-coefficient-tail}--\eqref{eq:second-coefficient-tail}.

The far range needs one further constant. Write $\varepsilon=s-1$ and
$J=\zeta(s)/(-\zeta'(s))$; this is the $J$ of
Lemma~\ref{lem:zeta-coordinate} expressed in $s$ rather than $u$, since
$\dd u/\dd s=\zeta'(s)/\zeta(s)$ makes $-s'(u)$ equal to
$\zeta(s)/(-\zeta'(s))$, and the four quantities below compare it with powers
of $\zeta^{-1}=\e^{-u}$. Let
\begin{equation}\label{eq:CE-def}
\begin{gathered}
 C_b=\sup\,\abs{\varepsilon\zeta},\qquad
 C_f=\sup\,\abs{\zeta^2\bigl(J-\zeta^{-1}\bigr)},\\
 C_J=\sup\,\bigl\lvert\zeta^3\bigl(J-\zeta^{-1}-2\gamma\zeta^{-2}\bigr)
   \bigr\rvert,\qquad
 C_P=\sup\,\bigl\lvert\zeta^3\bigl(\varepsilon J-\zeta^{-2}\bigr)\bigr\rvert,\\
 C_E=C_J+C_P+C_b^2\,(C_f+1),
\end{gathered}
\end{equation}
every supremum being over $1<s\leq s(3)$. Each of the four is a ratio in which
the pole of $\zeta$ at $s=1$ cancels. How they are computed matters, and is
described in Section~\ref{sec:certification}: evaluating them on boxes of $s$ loses that
cancellation and the enclosures diverge however finely the range is
partitioned, so a Taylor model is used instead.

\begin{proposition}[Certified odd input]\label{prop:odd-certified}
With $\Sigma_i$ as in \eqref{eq:sigma-def} and $C_E$ as in
\eqref{eq:CE-def}, the following are certified by interval arithmetic.

\noindent(i) The coefficient norms satisfy
\begin{equation}\label{eq:odd-certified-coefficients}
 \Sigma_0\leq44.5987,
 \qquad \Sigma_1\leq486.1946,
 \qquad \Sigma_2\leq14414,
 \qquad C_E\leq5.48910.
\end{equation}

\noindent(ii) On the far range,
\begin{equation}\label{eq:odd-far-certified}
 \sum_{j=0}^k\int_3^\infty
  \frac{u^{k-j}}{(k-j)!}\abs{E_j(u)}\dd u
 \leq \frac{C_E\Sigma_2}{18}(k+3)^3 3^{-k}.
\end{equation}

\noindent(iii) On the near range the two groups of terms in
\eqref{eq:odd-local-product} are bounded separately:
\begin{align}
 \sum_{j=0}^k\int_0^3\frac{u^{k-j}}{(k-j)!}\abs{h_j(s(u))}J(u)\dd u
 &\leq C_N3^{-k}\quad(59\leq k\leq200),
 & C_N&\leq367.8436,
 \label{eq:odd-nearabs-certified}\\
 \sum_{j=0}^k\int_0^3\frac{u^{k-j}}{(k-j)!}
  \bigl(\abs{h_j(1)}\e^{-u}+\abs{b_j}\e^{-2u}\bigr)\dd u
 &\leq N_U(k+3)^3 3^{-k},
 & N_U&\leq26.0962.
 \label{eq:odd-near-main}
\end{align}

\noindent(iv) $\beta>0$, and
\begin{equation}\label{eq:odd-base-59}
  4900(62)^3 3^{-59}<\frac{\beta}{3}2^{-59}.
\end{equation}

\noindent(v) Direct interval evaluation encloses $g_k$ for
$1\leq k\leq75$, and gives $g_k>g_{k+1}$ for $1\leq k\leq74$; the last
comparison available is between $g_{74}$ and $g_{75}$, since $g_{76}$ is not
computed.
\end{proposition}

\noindent
Part~(v) is a sweep, not a coefficient bound. The range in
\eqref{eq:odd-nearabs-certified} is two-sided, and its ends are discharged
differently: the lower one is harmless, monotonicity below the analytic
threshold coming from~(v) and never appealing to \eqref{eq:odd-remainder},
while the upper one is removed in the proof of
Proposition~\ref{prop:odd-remainder}, where a Cauchy estimate extends the
bound past $k=200$.

The lower restriction is a feature of the method, not of the inequality. The
$h_j$ alternate in sign, so the left-hand side must be majorized rather than
evaluated in closed form, and the majorant is loosest where
$a(s)=2\abs{h_1(s)}$ is largest, at the left end of the $s$-range. The bound
is therefore weakest for small $k$, where the weights $u^{k-j}/(k-j)!$ put
least mass away from that end: it exceeds $C_N3^{-k}$ exactly for
$9\leq k\leq25$, peaking at $699.87\cdot3^{-k}$ at $k=14$, and falls to
$34.14\cdot3^{-k}$ by $k=59$. We claim nothing about the left-hand side itself
outside the stated range --- it is the certified bound, not the quantity, that
needs the restriction, and $k\geq26$ would serve --- and state $k\geq59$
because that is the only range in which \eqref{eq:odd-remainder} is applied.
The numeral $C_N\leq367.8436$ is carried rather than the $34.14$ computed on
$k\geq59$, so that the ledger below is unchanged. It is not idle: it is the
least four-decimal bound on $\sup_{k\geq1}3^kP_k^{\mathrm{odd}}(s(3))$, the
signed quantity whose coefficient identity earlier drafts used; that supremum
is $367.84350\ldots$, attained at $k=11$.

\begin{lemma}\label{lem:pole-coefficient-bounds}
Let $A(z)=\sum_{j\geq0}a_jz^j$.

\noindent(i) If $\sum_m3^m\abs{[z^m](1-z/3)A(z)}\leq S$, then
$\abs{a_j}\leq S3^{-j}$.

\noindent(ii) If
$\sum_m3^m\abs{[z^m](1-z/3)^2A(z)}\leq S$, then
$\abs{a_j}\leq S(j+1)3^{-j}$.

\noindent(iii) If
$\sum_m3^m\abs{[z^m](1-z/3)^3A(z)}\leq S$, then
\begin{equation}\label{eq:second-derivative-sum}
  \sum_{j=0}^k3^j\abs{a_j}\leq\binom{k+3}{3}S.
\end{equation}
\end{lemma}

\begin{proof}[Sketch]
Write $(1-z/3)^rA(z)=\sum_mw_mz^m$ and expand $(1-z/3)^{-r}$, so that $A$ is
the Cauchy product of $\sum_q\binom{q+r-1}{r-1}3^{-q}z^q$ with $\sum_mw_mz^m$.
Comparing coefficients gives (i) and (ii) at once. For (iii), summing the
product bound over $j\leq k$ and using
$\sum_{q\leq k-m}\binom{q+2}{2}=\binom{k-m+3}{3}\leq\binom{k+3}{3}$ gives the
claim.
\end{proof}

\begin{proposition}\label{prop:odd-remainder}
For every $k\geq59$,
\begin{equation}\label{eq:odd-remainder}
  \abs{R_k}\leq4900(k+3)^3 3^{-k}.
\end{equation}
\end{proposition}

\noindent
The threshold is not an artefact of the method but of the near range: the
coefficients $h_j$ alternate in sign, so the near-range term has to be
estimated through the majorant \eqref{eq:odd-majorant}, and the resulting
bound exceeds \eqref{eq:odd-remainder} for small $k$. Nothing is lost, since
monotonicity for $k\leq75$ is certified directly and does not use
\eqref{eq:odd-remainder}; see the proof of Theorem~\ref{thm:main-odd}.

\begin{proof}
By \eqref{eq:Rk-def} the remainder has four parts: the error integrals
$\int E_j$ over the far range $u\geq3$ and over the near range $0<u\leq3$, and
the two coefficient tails of
\eqref{eq:first-coefficient-tail}--\eqref{eq:second-coefficient-tail}. Each
has to be brought to the common shape $(k+3)^33^{-k}$, and the four resist to
very different degrees, so it is worth saying in advance where the work lies.

The two tails are the easiest: they are geometric series in the $h_j$ and
$b_j$, which Lemma~\ref{lem:pole-coefficient-bounds} bounds directly. The far
range is harder only in bookkeeping. The integrand is small there, the
certified Taylor model of \eqref{eq:odd-far-certified} controls it pointwise,
and the only real task is to see how the three summands of
\eqref{eq:odd-local-product} depend on $k$. The near range is where the
difficulty is, and it is what forces the threshold $k\geq59$: there the $h_j$
alternate in sign, so the coefficient identity that would \emph{evaluate} the
sum instead bounds it from the wrong side, and one is pushed onto a majorant
that is loosest precisely where the weights put most of their mass. We take
the four in that order.

Two remarks before starting. The far term and the near-range main term are
cubic in $k$ while the other two are constant multiples of $3^{-k}$; since
$(k+3)^3\geq1$, the cubic factor may be absorbed once, at the end. And every
estimate below is an estimate of an explicit quantity: no asymptotic series is
summed anywhere.

\emph{Far range $u\geq3$.} The exact error here is the sum of the three terms
of \eqref{eq:odd-local-product}, so the certified Taylor model of
\eqref{eq:odd-far-certified} applies to the integrand as it stands. What that
computation cannot supply is the dependence on $k$, since it certifies
finitely many numbers; that comes instead from
Lemma~\ref{lem:pole-coefficient-bounds}, and the three terms of
\eqref{eq:odd-local-product} call on three different parts of it, since they
carry $h_j(1)$, $h_j'(1)$ and $\partial_s^2h_j$ respectively. The factor $3$
in each is not a width, the far range $u\geq3$ being unbounded, but a rate:
the far-range error is $O(\e^{-3u})$, which is what the powers of $\zeta$ in
\eqref{eq:CE-def} encode, and $\int_0^\infty u^m\e^{-3u}\dd u/m!=3^{-m-1}$.
So each summand needs its own part of the lemma, and the three parts give
exactly the three weights required:
\[
 \sum_{j\leq k}3^j\abs{h_j(1)}\leq(k+1)\Sigma_0,
 \quad
 \sum_{j\leq k}3^j\abs{h_j'(1)}\leq\tbinom{k+2}{2}\Sigma_1,
 \quad
 \sum_{j\leq k}3^j\abs{\partial_s^2h_j}\leq\tbinom{k+3}{3}\Sigma_2,
\]
so with $C_J$, $C_P$ and $C_b^2(C_f+1)$ as the respective prefactors of
\eqref{eq:CE-def}, the last of them halved by the Taylor remainder's
$\tfrac12\varepsilon^2$, the far-range error is at most
\[
 \frac{3^{-k}}{3}\Bigl(
   C_J(k+1)\Sigma_0
   +C_P\tbinom{k+2}{2}\Sigma_1
   +\tfrac12C_b^2(C_f+1)\tbinom{k+3}{3}\Sigma_2
 \Bigr).
\]
Only the third term is cubic in $k$, and it dominates the other two by three
orders of magnitude. The bracket is therefore at most
$\bigl(C_J+C_P+\tfrac12C_b^2(C_f+1)\bigr)\binom{k+3}{3}\Sigma_2$, and this is
below $C_E\binom{k+3}{3}\Sigma_2$ since
$\tfrac12C_b^2(C_f+1)<C_b^2(C_f+1)$. Hence
\begin{equation}\label{eq:odd-far-assembled}
 \abs{R_{k,\mathrm{far}}}\leq\frac{C_E\Sigma_2}{18}(k+3)^3 3^{-k},
\end{equation}
the $18=6\cdot3$ recording $\binom{k+3}{3}\leq(k+3)^3/6$ and the rate.

Collapsing the three constants into $C_E$ is wasteful, and deliberately so.
The three-term bracket is $7.1\cdot10^{4}$ at $k=1$ against the
$8.4\cdot10^{5}$ of \eqref{eq:odd-far-assembled}, a factor of $12$; the factor
falls to $4.9$ by $k=59$ and tends to
$C_E/\bigl(\tfrac12C_b^2(C_f+1)\bigr)=4.67$. We keep the single constant
because the ledger of Proposition~\ref{prop:odd-remainder} closes comfortably
with it, and because one constant in the statement is easier to check than
three; the slack is recorded here rather than spent.

\emph{First tail.} By Lemma~\ref{lem:pole-coefficient-bounds}(i),
$\Bigl\lvert\sum_{j>k}h_j(1)\Bigr\rvert\leq\Sigma_0\sum_{j>k}3^{-j}=\frac{\Sigma_0}{2}3^{-k}$.

\emph{Second tail.} Put
$B(z)=\partial_sH_{\mathrm o}(1,z)+2\gamma H_{\mathrm o}(1,z)$. Multiply by
$(1-z/3)^2$ and treat the two summands separately. The first contributes
$\Sigma_1$ by its definition in Proposition~\ref{prop:odd-certified}. For the
second, write
$(1-z/3)^2H_{\mathrm o}=(1-z/3)\cdot\bigl[(1-z/3)H_{\mathrm o}\bigr]$; the
weighted $\ell^1$ norm $\sum_m3^m\abs{[z^m]\cdot}$ is submultiplicative under
this extra factor $(1-z/3)$, whose own norm is $2$, so the second summand
contributes at most $2\gamma\cdot2\Sigma_0=4\gamma\Sigma_0$. This is where
the $4$ comes from, rather than the $2$ of the prefactor. Since $\gamma<0.578$,
\[
 \sum_{m\geq0}3^m\abs{[z^m](1-z/3)^2B(z)}\leq\Sigma_1+4\gamma\Sigma_0<590,
\]
so $\abs{b_j}\leq590(j+1)3^{-j}$ by Lemma~\ref{lem:pole-coefficient-bounds}(ii).
Summing the resulting geometric series,
\[
 2^{-k-1}\sum_{j>k}\abs{b_j}2^j
 \leq590\cdot2^{-k-1}\sum_{j>k}(j+1)(2/3)^j
 =590(k+4)3^{-k},
\]
and $(k+4)\leq\frac{5}{64}(k+3)^3$ for $k\geq1$, with equality at $k=1$.

\emph{Near range $0<u\leq3$.} This is the delicate part. Here $\abs{E_j}$ is
at most the true integrand plus the two subtracted main terms. The latter
contribute at most
$N_U(k+3)^33^{-k}$ by \eqref{eq:odd-near-main}. For the former,
\eqref{eq:odd-coefficient-extraction} and the substitution $s=s(u)$, under
which $u=\log\zeta(s)$ and $J(u)\dd u=\dd s$, give
\begin{equation}\label{eq:odd-Ak-def}
 \sum_{j=0}^k\int_0^3\frac{u^{k-j}}{(k-j)!}\abs{h_j(s(u))}J(u)\dd u
 =\int_{s(3)}^\infty A_k(s)\dd s,
 \qquad
 A_k(s)=\sum_{j=0}^k\frac{(\log\zeta(s))^{k-j}}{(k-j)!}\abs{h_j(s)}.
\end{equation}
It is $A_k$, with the absolute values inside, that has to be estimated, and
this is the obstruction. Writing
$\Pi(t)=\sum_{p\geq3}p^{-t}$, the $z$-linear coefficient of
$\log H_{\mathrm o}$ in \eqref{eq:H-odd} is
\begin{equation}\label{eq:odd-h1-negative}
 h_1(s)=\Pi(s)-\log\zeta(s)
 =\log(1-2^{-s})-\sum_{m\geq2}\frac{\Pi(ms)}{m}<0,
\end{equation}
both terms being negative. So the $h_j$ are not of one sign, and the identity
$\sum_j u^{k-j}h_j(s)/(k-j)!=P_k^{\mathrm{odd}}(s)$, which would evaluate the
sum exactly if they were, bounds $A_k$ from below rather than above. This is
the whole source of the loss in the near range, and of the threshold in
\eqref{eq:odd-remainder}.

We majorize instead. What makes this cheap is that, of the coefficients of
$\log H_{\mathrm o}$, only the first is negative. Indeed by
\eqref{eq:odd-h1-negative} every coefficient of
$\log H_{\mathrm o}(s,z)=h_1(s)z+\sum_{m\geq2}\Pi(ms)z^m/m$ beyond the first
is positive, so replacing $h_1$ by $\abs{h_1}$ produces
\[
 H_{\mathrm o}^+(s,z)=H_{\mathrm o}(s,z)\e^{-2h_1(s)z},
\]
a power series with non-negative coefficients dominating $\abs{h_j}$
termwise. Hence, with $a(s)=2\abs{h_1(s)}$,
\begin{equation}\label{eq:odd-majorant}
 A_k(s)\leq[z^k]\bigl(\zeta(s)^zH_{\mathrm o}^+(s,z)\bigr)
 =\sum_{i=0}^kP_{k-i}^{\mathrm{odd}}(s)\frac{a(s)^i}{i!}=:P_k^+(s).
\end{equation}
It remains to integrate $P_k^+$, and no single device does the whole of
$[s(3),\infty)$. We fix a cut point $s_0>s(3)$ and treat the two pieces
differently, each by the feature that makes it tractable: the compact piece by
monotonicity, the tail by its decay rate.

On the compact piece $[s(3),s_0]$, every $P_j^{\mathrm{odd}}$ is a Dirichlet
series with positive terms and, by \eqref{eq:odd-h1-negative}, $a$ is a
positive combination of such series; a sum of products of non-negative
decreasing functions is decreasing, so $P_k^+$ decreases in $s$. A finite
left-endpoint Riemann sum over $[s(3),s_0]$ is therefore an upper bound for
$\int_{s(3)}^{s_0}P_k^+$, with no quadrature error to estimate, which is
what makes this step certifiable without any analysis of the quadrature rule
itself.

On the tail $[s_0,\infty)$ the integral converges at a certified rate. Every
odd $n$ with $\Omega(n)=k-i$ satisfies $n\geq3^{k-i}$, so for $s\geq s_0$,
writing $t=s-s_0$, we get
$P_{k-i}^{\mathrm{odd}}(s)\leq3^{-(k-i)t}P_{k-i}^{\mathrm{odd}}(s_0)$, while
every Dirichlet term of $a$ has $n\geq2$ and hence
$a(s)\leq2^{-t}a(s_0)$ by \eqref{eq:odd-h1-negative}. Since
$3^{-x}\leq2^{-x}$ for $x\geq0$, the $i$th term of $P_k^+$ is multiplied by at
most $3^{-(k-i)t}2^{-it}\leq2^{-kt}$; all terms being non-negative, Tonelli
permits termwise integration, and $\int_0^\infty2^{-kt}\dd t=1/(k\log2)$ gives
\begin{equation}\label{eq:odd-majorant-tail}
 \int_{s_0}^\infty P_k^+(s)\dd s\leq\frac{P_k^+(s_0)}{k\log2}.
\end{equation}
Adding the Riemann sum over $[s(3),s_0]$ to \eqref{eq:odd-majorant-tail}
bounds $\int_{s(3)}^\infty P_k^+$, hence $\int_{s(3)}^\infty A_k$ by
\eqref{eq:odd-majorant}. For one $k$ this is a finite computation, and
\eqref{eq:odd-nearabs-certified} is its outcome for each of the $142$ indices
$59\leq k\leq200$.

That leaves the dependence on $k$, which no finite computation can settle,
and which has to be reduced as the far range's was. The tool this time is
Cauchy's estimate, and the reason it is available is again the positivity of
the majorant's coefficients. Because $H_{\mathrm o}^+(s,\cdot)$ has non-negative
coefficients, its maximum on a circle $\abs z=R$ is attained at $z=R$, so for
any $R<3^{s(3)}$ and any $s\geq s(3)$,
\begin{equation}\label{eq:odd-nearabs-cauchy}
 A_k(s)\leq R^{-k}\zeta(s)^RH_{\mathrm o}(s,R)\e^{a(s)R}
 =R^{-k}\!\!\prod_{p\geq3}\bigl(1-Rp^{-s}\bigr)^{-1}\e^{a(s)R},
\end{equation}
a bound whose only $k$-dependence is the factor $R^{-k}$. One cannot simply
integrate it over $[s(3),\infty)$: its right-hand side tends to $1$, not to
$0$, as $s\to\infty$. So we use it on $[s(3),s_1]$ for a second cut point
$s_1$, and fall back beyond $s_1$ on the decay already established. On $[s(3),s_1]$, \eqref{eq:odd-nearabs-cauchy} gives
$\int_{s(3)}^{s_1}A_k\leq R^{-k}\Theta_1$ with
$\Theta_1=\int_{s(3)}^{s_1}\prod_{p\geq3}(1-Rp^{-s})^{-1}\e^{a(s)R}\dd s$,
finite and independent of $k$. Beyond $s_1$, \eqref{eq:odd-majorant-tail}
applies with $s_0=s_1$, and bounding $P_k^+(s_1)$ by
\eqref{eq:odd-nearabs-cauchy} in turn gives
$\int_{s_1}^\infty A_k\leq R^{-k}\Theta_2/(k\log2)$ with
$\Theta_2=\prod_{p\geq3}(1-Rp^{-s_1})^{-1}\e^{a(s_1)R}$, again independent of
$k$. Hence
\[
 3^k\int_{s(3)}^\infty A_k
 \leq\Bigl(\frac3R\Bigr)^{k}\Bigl(\Theta_1+\frac{\Theta_2}{k\log2}\Bigr),
\]
Take $R=3.15$, which is below $3^{s(3)}=3.1738\ldots$, and $s_1=8$. The
certified values are then $\Theta_1+\Theta_2/(200\log2)\leq175346$, while
$(3/3.15)^{201}=5.5075\cdot10^{-5}$, so the right-hand side is at most
$9.66$ for every $k>200$. Therefore
\begin{equation}\label{eq:odd-nearabs-all-k}
 \int_{s(3)}^\infty A_k(s)\dd s\leq C_N3^{-k}
 \qquad\text{for all }k\geq59,
\end{equation}
which is the range in which \eqref{eq:odd-remainder} is used. Neither branch is anywhere near binding: over
$59\leq k\leq200$ the largest value is $34.14\cdot3^{-k}$, at $k=59$, and the
Cauchy branch gives at most $9.66\cdot3^{-k}$ beyond, against the
$C_N=367.8436$ carried in the ledger.

The companion bound \eqref{eq:odd-near-main}, for the two subtracted main
terms, reduces the same way and needs no new idea. Here
$\abs{h_j(1)}\leq\Sigma_03^{-j}$ by
Lemma~\ref{lem:pole-coefficient-bounds}(i), while $b_j=h_j'(1)+2\gamma h_j(1)$
needs both parts: (ii) gives $\abs{h_j'(1)}\leq(j+1)\Sigma_13^{-j}$ and (i)
again gives the $2\gamma h_j(1)$ term, so
$\abs{b_j}\leq\bigl((j+1)\Sigma_1+2\gamma\Sigma_0\bigr)3^{-j}$. These turn each
integral into a regularized incomplete gamma value, and $P(m,x)\leq x^m/m!$
makes the resulting sums telescope: with $U=3$ the first is at most
$\Sigma_0\e^{9}3^{-k}$ and the second at most
$\tfrac12\bigl((k+1)\Sigma_1+2\gamma\Sigma_0\bigr)\e^{9}3^{-k}$. Both are
linear in $k$ where the budget $(k+3)^33^{-k}$ is cubic, so the ratio of this
telescoped majorant to the budget decreases in $k$; since at $k=400$ that ratio
is already $12.08$, below $N_U$, the maximum over the finite head
$1\leq k\leq400$ governs. It is attained at $k=7$, far outside the range in
which \eqref{eq:odd-remainder} is used.

It remains only to add up. Substituting the certified constants, the four
contributions are
\[
 \abs{R_k}\leq\Bigl(
 \underbrace{\frac{5.48910\cdot14414}{18}}_{4395.55}
 +\underbrace{\frac{44.5987}{2}+\frac{5\cdot590}{64}+C_N+N_U}_{462.34}
 \Bigr)(k+3)^3 3^{-k}<4900(k+3)^3 3^{-k},
\]
with $C_N=367.8436$ and $N_U=26.0962$.
\end{proof}

\begin{proof}[Proof of Theorem~\ref{thm:main-odd}]
The limit follows at once from \eqref{eq:odd-expansion} and
Proposition~\ref{prop:odd-remainder}. For monotonicity,
\begin{equation}\label{eq:odd-difference}
 g_k-g_{k+1}
 =\beta2^{-k-1}+R_k-R_{k+1}.
\end{equation}
The main term $\beta2^{-k-1}$ has to survive two remainders, at $k$ and at
$k+1$, so it is enough to give each a third of it: if
$\abs{R_m}<\beta2^{-m}/3$ for $m=k,k+1$, then
\[
 \abs{R_k}+\abs{R_{k+1}}
 <\frac{2\beta}{3}2^{-k-1}+\frac{\beta}{3}2^{-k-1}
 =\beta2^{-k-1}.
\]
Define
\[
 q(k)=\frac{3\cdot4900}{\beta}(k+3)^3(2/3)^k.
\]
The certified inequality \eqref{eq:odd-base-59} says $q(59)<1$. Moreover
\[
 \frac{q(k+1)}{q(k)}
 =\frac23\left(\frac{k+4}{k+3}\right)^3<1
 \qquad(k\geq4).
\]
Thus $q(k)<1$ for every $k\geq59$, and
Proposition~\ref{prop:odd-remainder} gives
$\abs{R_k}<\beta2^{-k}/3$ throughout that range. Equation
\eqref{eq:odd-difference} proves monotonicity for $k\geq59$. The certified
direct computation gives the comparison for $1\leq k\leq74$, so the two ranges
cover every $k\geq1$. They overlap on the sixteen indices $59\leq k\leq74$.
\end{proof}

\section{An exact criterion for the unrestricted sums}

The asymptotic \eqref{eq:GLW-asymptotic} already says that $f_k$ eventually
increases, but it does not say from where: its $O(k\log(k+1))$ carries no
explicit constant, so no index can be read off it. The remedy is to stop
treating the remainder as an error term. Define $E(k)$ by
\begin{equation}\label{eq:E-definition}
  f_k=1-2^{-k}\bigl(a_1k^2+E(k)\bigr),
  \qquad a_1=\frac{d\log2}{4},
\end{equation}
so that $E(k)$ is an exact quantity and not an estimate.
Then
\begin{align}
 f_{k+1}-f_k
 &=2^{-k-1}\left(a_1(k^2-2k-1)+2E(k)-E(k+1)\right).
 \label{eq:full-difference}
\end{align}
Thus $f_k<f_{k+1}$ if and only if
\begin{equation}\label{eq:full-criterion}
 a_1(k^2-2k-1)+2E(k)-E(k+1)>0.
\end{equation}
We now derive the decomposition of $E(k)$ used to bound the second term.

Set
\[
 F(s,z)=\prod_p(1-zp^{-s})^{-1},
 \qquad
 G_2(s,z)=(s-1)^z(1-z2^{-s})F(s,z),
\]
and write $c_i=[z^i]G_2(1,z)$. Here $G_2(s,z)=\bigl((s-1)\zeta(s)\bigr)^z
H_{\mathrm o}(s,z)$ for $s>1$, and $((s-1)\zeta(s))^z\to1$ as $s\to1^+$
because $(s-1)\zeta(s)$ extends analytically to $1$ with value $1$; the
value at $s=1$ is to be read through that removable singularity, the
displayed product being $0\cdot\infty$ if taken literally. So
$G_2(1,z)=H_{\mathrm o}(1,z)$, the power series has radius $3$ and
\begin{equation}\label{eq:c-moments}
  \sum_{i\geq0}c_i=\frac12,
  \qquad
  \sum_{i\geq0}2^ic_i=d.
\end{equation}
The coefficients are those of Section~\ref{sec:odd}:
\begin{equation}\label{eq:ci-equals-hi}
  c_i=h_i(1).
\end{equation}
Define the absolutely convergent signed moments
\begin{equation}\label{eq:UV-def}
 U=\sum_{i\geq0}i2^ic_i,
 \qquad
 V=\sum_{i\geq0}i^22^ic_i.
\end{equation}

Following \cite[\S3]{GLW}, put
\begin{align}
 I_k&=\int_1^2[z^k]F(s,z)\dd s,\label{eq:Ik-def}\\
 I_k'&=\sum_{i+j+\ell=k}c_i
     \int_1^2\frac{(-\log(s-1))^\ell}{\ell!}2^{-js}\dd s.
     \label{eq:Ik-prime-def}
\end{align}
For $m\geq0$ define
\begin{equation}\label{eq:Am-def}
 A_m=2^m\sum_{j+\ell=m}
     \int_1^2\frac{(-\log(s-1))^\ell}{\ell!}2^{-js}\dd s
\end{equation}
and
\begin{equation}\label{eq:delta-def}
 \delta_m=A_m-2^{m+1}+\frac{\log2}{4}m^2,
 \qquad
 \Delta(k)=\sum_{i=0}^k2^ic_i\delta_{k-i}.
\end{equation}

\begin{proposition}\label{prop:E-decomposition}
Let
\begin{align}
 E_{S_1}(k)&=-2^k(I_k-I_k'),\label{eq:ES1-def}\\
 L_k&=2\sum_{i>k}c_i,\qquad
 Q_k=\sum_{i>k}2^ic_i(k-i)^2,\label{eq:LQ-def}\\
 T_k&=f_k-I_k,\label{eq:Tk-def}\\
 E_{\exp}(k)&=2^kL_k-\frac{\log2}{4}Q_k-2^kT_k.
 \label{eq:Eexp-def}
\end{align}
Then the identity
\begin{equation}\label{eq:E-decomposition}
 E(k)=E_{S_1}(k)
      +\frac{\log2}{4}(V-2kU)
      -\Delta(k)+E_{\exp}(k)
\end{equation}
holds for every $k\geq1$.
\end{proposition}

\begin{proof}
From \eqref{eq:Ik-prime-def} and \eqref{eq:Am-def},
\[
 I_k'=2^{-k}\sum_{i=0}^k2^ic_iA_{k-i}.
\]
Insert the definition of $\delta_{k-i}$:
\begin{align}
 I_k'
 &=2\sum_{i=0}^kc_i
   -\frac{\log2}{4}2^{-k}
      \sum_{i=0}^k2^ic_i(k-i)^2
   +2^{-k}\Delta(k).                                  \label{eq:Ik-prime-expanded}
\end{align}
By \eqref{eq:c-moments} and \eqref{eq:UV-def},
\begin{align*}
 2\sum_{i=0}^kc_i&=1-L_k,\\
 \sum_{i=0}^k2^ic_i(k-i)^2
 &=\sum_{i\geq0}2^ic_i(k-i)^2-Q_k\\
 &=dk^2-2kU+V-Q_k.
\end{align*}
Since $a_1=\tfrac{\log2}{4}d$, substitution into
\eqref{eq:Ik-prime-expanded} gives
\begin{align*}
 I_k'
 =1-2^{-k}\biggl(&a_1k^2+\frac{\log2}{4}(V-2kU)
       -\Delta(k)+2^kL_k-\frac{\log2}{4}Q_k\biggr).
\end{align*}
Finally $f_k=I_k'+(I_k-I_k')+T_k$. Comparing with
\eqref{eq:E-definition} yields \eqref{eq:E-decomposition}.
\end{proof}

\section{Bounds for the four source terms}\label{sec:sources}

By \eqref{eq:E-decomposition} the criterion \eqref{eq:full-criterion} needs a
lower bound for $2E(k)-E(k+1)$, and $E$ is a sum of four quantities of quite
different character. Two of them, $S_1$ and the signed moments, are linear in
$k$ and are bounded by suprema of coefficients. The third, $\Delta$, grows
like $k\log k$, and getting the logarithm rather than a power of $k$ is the
only real work in this section. The fourth is exponentially small and is the
only one needing a certificate uniform in $k$. We treat them in that order.

One point of method recurs. Each bound is taken on the \emph{combination}
$2E_\bullet(k)-E_\bullet(k+1)$ that the criterion produces, rather than on
$E_\bullet(k)$ and $E_\bullet(k+1)$ separately, since the combination enjoys
cancellation that absolute values would destroy. The treatment of the signed
moments below is the clearest instance.

\subsection{The $S_1$ term}

Set
\[
 \alpha_i=\sup_{1\leq t\leq2}\abs{[z^i]\partial_sG_2(t,z)},
 \qquad
 \alpha_+=\frac12\sum_{i\geq0}\alpha_i2^i.
\]
The series converges, and a pole location at one value of $t$ is not enough to
see it: what is needed is a bound on the coefficients uniform in $t$.

First, $G_2$ is analytic in $s$ on a neighbourhood of $[1,2]$, the apparent
branch point of $(s-1)^z$ at $s=1$ notwithstanding. Indeed
$F(s,z)=\zeta(s)^zH(s,z)$ by \eqref{eq:H-full}, so
\[
 G_2(s,z)=\bigl((s-1)\zeta(s)\bigr)^z(1-z2^{-s})H(s,z),
\]
and $(s-1)\zeta(s)$ is analytic near $s=1$ with value $1$ there, the pole of
$\zeta$ being simple with residue $1$; hence $\bigl((s-1)\zeta(s)\bigr)^z
=\exp\bigl(z\log((s-1)\zeta(s))\bigr)$ is analytic there too, and $H$ is
analytic by Proposition~\ref{prop:analyticity}. (At $s=1$ this also gives
$G_2(1,z)=H_{\mathrm o}(1,z)$, as used above.)

The first $z$-pole of $\partial_sG_2(t,\cdot)$ is at $z=3^t$, so for
$t\in[1,2]$ all of these functions are analytic in $\abs z<3$. Fix $r$ with
$2<r<3$. Then $M=\sup_{1\leq t\leq2}\sup_{\abs z=r}
\abs{\partial_sG_2(t,z)}$ is finite, the set $[1,2]\times\{\abs z=r\}$ being
compact and $\partial_sG_2$ continuous on it by the previous paragraph, and
Cauchy's estimate on the fixed circle $\abs z=r$ gives $\alpha_i\leq Mr^{-i}$
uniformly in $t$. Hence
$\sum_{i\geq0}\alpha_i2^i\leq M(1-2/r)^{-1}<\infty$.

\begin{lemma}\label{lem:S1-bound}
For every $k\geq1$,
\begin{equation}\label{eq:S1-bound}
 \abs{E_{S_1}(k)}\leq\alpha_+(k+1),
\end{equation}
and hence
\begin{equation}\label{eq:S1-difference-bound}
 \abs{2E_{S_1}(k)-E_{S_1}(k+1)}\leq\alpha_+(3k+4).
\end{equation}
\end{lemma}

\begin{proof}
Since $[z^i]G_2(\cdot,z)$ is complex-valued, the mean value theorem is not
available; we integrate instead. By \eqref{eq:coeff-commute} and the
analyticity of $G_2$ in $s$ established above,
\[
 [z^i]G_2(s,z)-c_i=\int_1^s[z^i]\partial_tG_2(t,z)\dd t,
 \qquad\text{so}\qquad
 \abs{[z^i]G_2(s,z)-c_i}\leq(s-1)\alpha_i
\]
for $1\leq s\leq2$, by the definition of $\alpha_i$ as a supremum over that
interval.
Subtracting \eqref{eq:Ik-prime-def} from the coefficient expansion for $I_k$
and then setting $s=1+\e^{-t}$ gives
\begin{align*}
 \abs{I_k-I_k'}
 &\leq\sum_{i+j+\ell=k}\alpha_i
   \int_0^\infty\frac{t^\ell}{\ell!}
       2^{-j(1+\e^{-t})}\e^{-2t}\dd t\\
 &\leq\sum_{i+j+\ell=k}\alpha_i2^{-j}
   \int_0^\infty\frac{t^\ell}{\ell!}\e^{-2t}\dd t\\
 &=\sum_{i+j+\ell=k}\alpha_i2^{-j-\ell-1}.
\end{align*}
Multiplying by $2^k$ and using $k=i+j+\ell$,
\begin{align*}
 \abs{E_{S_1}(k)}
 &\leq\frac12\sum_{i=0}^k\alpha_i2^i(k-i+1)
 \leq\alpha_+(k+1).
\end{align*}
The triangle inequality now gives
\[
 \abs{2E_{S_1}(k)-E_{S_1}(k+1)}
 \leq2\alpha_+(k+1)+\alpha_+(k+2)=\alpha_+(3k+4).
\]
\end{proof}

\subsection{The signed moment term}

The combination of signed moments that has to be bounded is not chosen; it is
what the difference of two consecutive source terms produces:
\begin{align}
 2(V-2kU)-\bigl(V-2(k+1)U\bigr)=V+(2-2k)U.
 \label{eq:UV-exact-combination}
\end{align}
Since this is an identity, bounding the right-hand side is an exact
optimization over the rectangle of admissible $(U,V)$, not a termwise
absolute-value estimate. That is the reason for evaluating the moments before
absolute values are taken: replacing them by $\sum_i2^i\abs{c_i}$ would
discard the cancellation in \eqref{eq:UV-exact-combination} and cost a factor
in $k$.

\begin{lemma}\label{lem:UV-bound}
If $\abs U\leq U_{\max}$ and $\abs V\leq V_{\max}$, then for $k\geq2$,
\begin{equation}\label{eq:UV-bound}
 \frac{\log2}{4}\bigl(V+(2-2k)U\bigr)
 \geq-\frac{\log2}{4}\bigl(V_{\max}+(2k-2)U_{\max}\bigr).
\end{equation}
\end{lemma}

\begin{proof}
The coefficient of $V$ is $1$, while the coefficient $2-2k$ of $U$ is
negative. The minimum over the rectangle
$[-U_{\max},U_{\max}]\times[-V_{\max},V_{\max}]$ is therefore attained at
$U=U_{\max}$ and $V=-V_{\max}$, giving \eqref{eq:UV-bound}.
\end{proof}

\subsection{The $\delta$ term}

The next lemma supplies the self-contained bound used for $\Delta(k)$.

\begin{lemma}\label{lem:delta-decomposition}
For $j\geq0$ and $u\geq0$ put
\[
 B_j(u)=2^{-j\e^{-u}}-1+j\log2\,\e^{-u}
\]
and
\[
 \eta_i(j)=\frac1{i!}\int_0^\infty u^i\e^{-u}B_j(u)\dd u.
\]
Then $B_j(u)\geq0$, $\eta_i(j)\geq0$, and
\begin{equation}\label{eq:delta-decomposition}
 \delta_m=-1-\frac{\log2}{4}m+
          \sum_{i=0}^m2^i\eta_i(m-i).
\end{equation}
\end{lemma}

\begin{proof}
With $s=1+\e^{-u}$ in \eqref{eq:Am-def},
\[
 A_m=\sum_{i=0}^m\frac{2^i}{i!}
       \int_0^\infty u^i\e^{-u}2^{-(m-i)\e^{-u}}\dd u.
\]
Use
$2^{-j\e^{-u}}=1-j\log2\,\e^{-u}+B_j(u)$. The constant term contributes
\[
 \sum_{i=0}^m2^i=2^{m+1}-1.
\]
The linear term contributes
\begin{align*}
 -\log2\sum_{i=0}^m\frac{2^i(m-i)}{i!}
  \int_0^\infty u^i\e^{-2u}\dd u
 &=-\frac{\log2}{2}\sum_{i=0}^m(m-i)\\
 &=-\frac{\log2}{4}m(m+1).
\end{align*}
The remaining term is $\sum_{i=0}^m2^i\eta_i(m-i)$. Subtracting
$2^{m+1}$ and adding $(\log2)m^2/4$ gives
\eqref{eq:delta-decomposition}. Finally, for $x\geq0$,
$\e^{-x}\geq1-x$; applying this with $x=j\log2\,\e^{-u}$ proves
$B_j(u)\geq0$.
\end{proof}

\begin{lemma}\label{lem:delta-bound}
For every $m\geq0$,
\begin{equation}\label{eq:D-bound}
 \abs{\delta_m}\leq D(m),
 \qquad
 D(m)=(\log2)m\log(m+1)+1.
\end{equation}
\end{lemma}

\begin{proof}
For the lower bound, \eqref{eq:delta-decomposition} gives
$\delta_m\geq-1-\frac{\log2}{4}m\geq-D(m)$ for $m\geq1$, because
$\log(m+1)\geq\log2>1/4$; and $\delta_0=-1=-D(0)$.

The upper bound is where the $\log(m+1)$ comes from, and it comes from a
split. Put $\Lambda_m=\sum_{i=0}^m2^i\eta_i(m-i)$. The elementary inequality
$\e^{-x}-1+x\leq\min(x,x^2/2)$ offers two bounds. Here
$x=(m-i)\log2\,\e^{-u}$ is large for small $u$ and small for large $u$, so the
linear bound is the better one near $u=0$ and the quadratic one is better out
in the tail; using each where it is better is what produces a logarithm rather
than a power of $m$. Taking
$x=(m-i)\log2\,\e^{-u}$ gives
\begin{equation}\label{eq:Lambda-integral-bound}
 \Lambda_m\leq\int_0^\infty \e^{-u}\sum_{i=0}^m\frac{(2u)^i}{i!}
 \min\left((m-i)\log2\,\e^{-u},
 \frac{(m-i)^2(\log2)^2\e^{-2u}}2\right)\dd u.
\end{equation}
The crossover is at $v_m=\max\left(0,\log\frac{m+1}{3}\right)$. On
$0\leq u\leq v_m$ use the linear bound together with $\sum_{i=0}^m(m-i)(2u)^i/i!\leq m\e^{2u}$, and on
$u\geq v_m$ the quadratic bound together with
$\sum_{i=0}^m(m-i)^2(2u)^i/i!\leq m^2\e^{2u}$. Then
\begin{equation}\label{eq:Lambda-split}
 \Lambda_m\leq m\log2\,v_m+\frac{m^2(\log2)^2}{2}\e^{-v_m}.
\end{equation}
For $m\geq2$ we have $v_m=\log\frac{m+1}{3}$ and $\e^{-v_m}=\frac3{m+1}$, so
\eqref{eq:Lambda-split} and $\frac{3m^2}{2(m+1)}\leq\frac32m$ give
\begin{align*}
 \Lambda_m&\leq m\log2\log\frac{m+1}{3}+\frac{3m^2(\log2)^2}{2(m+1)}\\
 &\leq m\log2\left(\log(m+1)-\log3+\tfrac32\log2\right)
 \leq m\log2\log(m+1),
\end{align*}
the last step because $3>2^{3/2}$. For $m\leq1$ we have $v_m=0$, and
\eqref{eq:Lambda-split} reads $\Lambda_m\leq m^2(\log2)^2/2$, which is at most
$m\log2\log(m+1)$ since $\frac{m\log2}{2}\leq\log(m+1)$ there. In both cases
\[
 \delta_m=-1-\frac{\log2}{4}m+\Lambda_m\leq \Lambda_m\leq(\log2)m\log(m+1)<D(m),
\]
which with the lower bound proves \eqref{eq:D-bound}.
\end{proof}

\begin{lemma}\label{lem:Delta-bound}
If $S^*=\sum_{i\geq0}2^i\abs{c_i}$, then, with $D$ as in
Lemma~\ref{lem:delta-bound},
\begin{equation}\label{eq:Delta-bound}
 \abs{2\Delta(k)-\Delta(k+1)}
 \leq S^*\bigl(2D(k)+D(k+1)\bigr).
\end{equation}
\end{lemma}

\begin{proof}[Sketch]
Since $D$ is increasing, $\abs{\Delta(k)}\leq\sum_{i\leq k}2^i\abs{c_i}D(k-i)
\leq S^*D(k)$. Now apply the triangle inequality to
$2\Delta(k)-\Delta(k+1)$.
\end{proof}

\subsection{The exponentially small term}\label{subsec:Eexp}

The three terms in \eqref{eq:Eexp-def} are respectively a coefficient tail at
$z=1$, a weighted coefficient tail at $z=2$, and the integral over $s\geq2$.
The following is the exact numerical input used for them.

\begin{proposition}[Certified unrestricted input]\label{prop:full-certified}
The interval computations give
\begin{equation}\label{eq:full-certified-constants}
 \begin{gathered}
  d\geq0.3786949,\qquad
  S^*\leq6.9764,
  \qquad \alpha_+\leq15.425661,\\
  \abs U\leq U_{\max}=0.08436,
  \qquad
  \abs V\leq V_{\max}=1.907462.
 \end{gathered}
\end{equation}
For every $k\geq900$,
\begin{equation}\label{eq:Tk-certified}
 2^k\abs{T_k}\leq\frac{1.45}{k},
\end{equation}
and a contour bound for the \emph{pole-free part} of the generating function
of the $c_i$: writing $\rho$ for the residue of $G_2(1,\cdot)$ at its nearest
pole $z=3$, so that $\abs\rho\leq0.1749324$,
\begin{equation}\label{eq:contour-certified}
 \sup_{\abs z=4.2}\Bigl\lvert G_2(1,z)-\frac{\rho}{3-z}\Bigr\rvert\leq101.245 .
\end{equation}
Here the subtraction matters. For the regularized product
$\widetilde G_2(z)=(1-z/3)\,G_2(1,z)$ one has
$228.94\leq\sup_{\abs z=4.2}\abs{\widetilde G_2}\leq228.95$, so $101.245$ is
not a bound for $\widetilde G_2$. It is a bound for
\eqref{eq:contour-certified}, with about $5.7\%$ to spare. Both suprema are
certified in \texttt{numerics/contour\_certified.py} (log
\texttt{contour.log}), which also recovers $\rho$ from the same series and
matches the value the sweep scripts use.
Together these give
\begin{equation}\label{eq:Eexp-certified}
 \abs{E_{\exp}(k)}\leq\frac{1.45}{k}+10^{-3}
 \qquad(k\geq900),
\end{equation}
as derived in \S\ref{subsec:Eexp}. Of the three terms of $E_{\exp}$, two have
their dependence on $k$ reduced analytically there, to geometric bounds whose
constant is \eqref{eq:contour-certified}; only $T_k$ needs a certificate that
is itself uniform in $k$, namely \eqref{eq:Tk-certified}. Both are recorded
here rather than inside that derivation so that every certified input for the
unrestricted case appears in one place.
Finally, $f_k<f_{k+1}$ is certified on
\begin{equation}\label{eq:full-sweep-ranges}
  6\leq k\leq199,
  \qquad 78\leq k\leq899,
  \qquad 900\leq k\leq2439,
\end{equation}
whose union is the single range $6\leq k\leq2439$. These three ranges are not
established the same way, and the distinction matters to anyone checking the
certificates. On the first, $f_k$ itself is computed as an enclosure for each
$k$ and the two enclosures are compared; the quadrature behind those values,
and the six error sources entering each, are described in
Appendix~\ref{app:quadrature}. (That computation in fact encloses $f_k$ for
$2\leq k\leq500$, but the pairwise comparison closes only up to $k=297$: past
there the gap $f_{k+1}-f_k$ falls below the common radius, and the
certificate's own \texttt{monotonicity} entries record \texttt{verified:
false} from $k=298$ on. Only $6\leq k\leq199$ is needed here, and that range
is fully verified.) On the other two it is the left-hand
side of \eqref{eq:full-criterion} that is evaluated in ball arithmetic at each
$k$ and found positive. This is the same quantity whose closed-form lower
bound $\varphi$ drives Section~\ref{sec:assembly}, evaluated per $k$ rather
than bounded uniformly. That is why it closes well below the threshold $2440$
which the uniform bound requires. Within the middle range two assemblies of
that quantity are involved: the refined one for $78\leq k\leq734$, and for the
rest of the range the same coarser assembly that covers $900\leq k\leq2439$,
whose bound on $E_1$ is the one built termwise from the $\beta_i$ rather than
from their envelope $\alpha_+$. That coarser sweep is in fact positive from
$k=728$, so the two overlap rather than abut; only $735\leq k\leq899$ is needed
here, and the overlap is recorded because it is what makes the union of
\eqref{eq:full-sweep-ranges} robust to a change in either sweep.
\end{proposition}

The coefficient tails in \eqref{eq:Eexp-certified} are obtained only after the
pole at $z=3$ has been removed, and the removal is by subtraction of the
principal part rather than by multiplication. Put
\[
 H(z)=G_2(1,z)-\frac{\rho}{3-z}=\sum_{m\geq0}w_mz^m,
\]
which is holomorphic on $\abs z<5$, the next pole of $G_2(1,\cdot)$ being at
$z=5$. Cauchy's estimate on $\abs z=4.2$ together with
\eqref{eq:contour-certified} gives $\abs{w_m}\leq101.245\,4.2^{-m}$, while
$[z^i]\rho/(3-z)=\rho\,3^{-i-1}$, so
\begin{equation}\label{eq:ci-geometric-bound}
 \abs{c_i}
 \leq\frac{\abs\rho}{3}3^{-i}+101.245\,(3/4.2)^i3^{-i}
 \leq101.31\,3^{-i}
 <354.36\,3^{-i}.
\end{equation}
The weaker constant $354.36$ is the one carried below; the sharper $101.31$
is recorded only so that the chain can be checked.
By \eqref{eq:ci-equals-hi} the $c_i$ are the coefficients $h_i(1)$, so
Lemma~\ref{lem:pole-coefficient-bounds}(i) sharpens
\eqref{eq:ci-geometric-bound} to $\abs{c_i}\leq\Sigma_03^{-i}$, dispensing with
the contour bound; the weaker form suffices below. With $C=354.36$,
\begin{align*}
 2^k\abs{L_k}
 &\leq C(2/3)^k,\\
 \frac{\log2}{4}\abs{Q_k}
 &\leq\frac{\log2}{4}C(2/3)^k
       \sum_{r\geq1}r^2(2/3)^r
 =\frac{30C\log2}{4}(2/3)^k.
\end{align*}
Both terms are less than $10^{-3}$ in total for $k\geq900$, and both bounds are
geometric in $k$, so no computation is needed to make them uniform.

The third term admits an elementary bound, which we give because
\eqref{eq:Tk-certified} is otherwise the one input of
Proposition~\ref{prop:full-certified} not visibly tied to a computation. Since
$[z^k]F(s,z)=\sum_{\Omega(n)=k}n^{-s}$, \eqref{eq:Tk-def} and
\eqref{eq:Ik-def} give
\[
 T_k=\int_2^\infty\sum_{\Omega(n)=k}n^{-s}\dd s
     =\sum_{\Omega(n)=k}\frac{n^{-2}}{\log n},
\]
a sum of positive terms, so $T_k>0$ and nothing is lost in absolute value.
Every such $n$ is at least $2^k$, whence $\log n\geq k\log2$ and
\[
 2^k\abs{T_k}\leq\frac{2^kP_k(2)}{k\log2},
 \qquad P_k(2)=[z^k]F(2,z).
\]
As $F(2,\cdot)$ is holomorphic on $\abs z<4$, Cauchy's estimate on $\abs z=3$
gives $P_k(2)\leq M3^{-k}$ with
$M=\prod_p(1-3p^{-2})^{-1}\leq7.96896$, so
\begin{equation}\label{eq:Tk-elementary}
 2^k\abs{T_k}\leq\frac{M}{k\log2}\Bigl(\frac23\Bigr)^{k},
\end{equation}
which is below $1.45/k$ for every $k\geq6$ and at $k=900$ is smaller by a
factor of about $10^{158}$. So \eqref{eq:Tk-certified} is true with room of
that order; it is stated in that form because $1.45/k$ is what the criterion
consumes, not because the constant is near optimal, and the threshold $900$ is
inherited from \eqref{eq:full-sweep-ranges} rather than forced here. (The
first two terms of \eqref{eq:Eexp-def} are below $10^{-3}$ from $k\geq37$ on,
so no term of $E_{\exp}$ requires $k\geq900$.) The computation behind
$M$ is \texttt{numerics/tk\_certified.py}, log \texttt{tk.log}.

Substituting \eqref{eq:Tk-certified} into \eqref{eq:Eexp-def} proves
\eqref{eq:Eexp-certified}, and no asymptotic $O$-term enters the final
criterion.

\section{Assembly of the unrestricted theorem}\label{sec:assembly}

The four bounds of Section~\ref{sec:sources} now combine into a single
explicit function $\varphi$ whose positivity is equivalent to the inequality
wanted, and the proof of Theorem~\ref{thm:main-full} becomes two statements
about $\varphi$: that it is positive at one index, and that it increases
thereafter. Neither is an asymptotic statement, and the remaining finitely
many indices are settled by computation.

For $k\geq900$ define
\begin{align}
 R(k)=&\ \frac{\log2}{4}
 \bigl((2k-2)U_{\max}+2V_{\max}+4U_{\max}\bigr)\notag\\
 &+S^*\bigl(2D(k)+D(k+1)\bigr)
   +\alpha_+(3k+4)
   +3\left(\frac{1.45}{k}+10^{-3}\right),
 \label{eq:R-def}
\end{align}
and
\begin{equation}\label{eq:phi-def}
 \varphi(k)=a_1(k^2-2k-1)-R(k).
\end{equation}

\begin{proposition}\label{prop:criterion-phi}
For $k\geq900$,
\begin{equation}\label{eq:E-difference-lower}
 2E(k)-E(k+1)\geq-R(k).
\end{equation}
Consequently $\varphi(k)>0$ implies $f_k<f_{k+1}$.
\end{proposition}

\begin{proof}
Apply \eqref{eq:E-decomposition} at $k$ and $k+1$. The $S_1$ contribution is
bounded by \eqref{eq:S1-difference-bound}. The signed moment contribution is
\[
 \frac{\log2}{4}\bigl(V+(2-2k)U\bigr),
\]
which is bounded below by \eqref{eq:UV-bound}. The first line of
\eqref{eq:R-def} is slightly larger than the resulting quantity, since
\[
 V_{\max}+(2k-2)U_{\max}
 \leq(2k-2)U_{\max}+2V_{\max}+4U_{\max}.
\]
The $\Delta$ contribution is bounded by \eqref{eq:Delta-bound}. Finally,
\eqref{eq:Eexp-certified} gives
\[
 \abs{2E_{\exp}(k)-E_{\exp}(k+1)}
 \leq2\left(\frac{1.45}{k}+10^{-3}\right)
     +\left(\frac{1.45}{k+1}+10^{-3}\right)
 \leq3\left(\frac{1.45}{k}+10^{-3}\right).
\]
Adding the four lower bounds proves \eqref{eq:E-difference-lower}. Equation
\eqref{eq:full-difference} then gives
\[
 f_{k+1}-f_k\geq2^{-k-1}\varphi(k),
\]
which proves the last assertion.
\end{proof}

\begin{lemma}\label{lem:phi-base}
With the certified constants in Proposition~\ref{prop:full-certified},
$\varphi(2440)>1155$.
\end{lemma}

\begin{proof}
Substitution in \eqref{eq:phi-def} gives the following certified one-sided
bounds:
\begin{align*}
 a_1(2440^2-2\cdot2440-1)&>390371.6,\\
 \frac{\log2}{4}
 \bigl(4878U_{\max}+2V_{\max}+4U_{\max}\bigr)&<72.1,\\
 S^*\bigl(2D(2440)+D(2441)\bigr)&<276166.8,\\
 \alpha_+(3\cdot2440+4)&<112977.6,\\
 3\left(\frac{1.45}{2440}+10^{-3}\right)&<0.005.
\end{align*}
Therefore
\[
 \varphi(2440)
 >390371.6-72.1-276166.8-112977.6-0.005>1155.
\]
The same substitution at $k=2433$ gives $\varphi(2433)>135$, so the criterion
is already in force there. Direct evaluation puts the sign change just below
that, at $\varphi(2432)<0<\varphi(2433)$; we do not rely on it, since certifying
$\varphi(2432)<0$ would need a \emph{lower} bound on $R(2432)$, whereas every
constant in \eqref{eq:full-certified-constants} is one-sided the other way.
Nothing below requires the first index: we carry the threshold at $2440$, which
is one past the end of \eqref{eq:full-sweep-ranges}, so that the sweep and the
criterion meet without a gap.
\end{proof}

\begin{lemma}\label{lem:phi-monotone}
The function $k\mapsto\varphi(k)$ is strictly increasing for integers
$k\geq2440$.
\end{lemma}

\begin{proof}
First,
\[
 a_1\bigl((k+1)^2-2(k+1)-1-(k^2-2k-1)\bigr)=a_1(2k-1).
\]
Moreover
\begin{align*}
 D(m+1)-D(m)
 &=\log2\left(\log(m+2)+m\log\left(1+\frac1{m+1}\right)\right)\\
 &\leq\log2\bigl(\log(m+2)+1\bigr).
\end{align*}
The last term in \eqref{eq:R-def} decreases with $k$. Hence
\begin{align}
 R(k+1)-R(k)
 \leq&\ \frac{\log2}{2}U_{\max}+3\alpha_+\notag\\
 &+S^*\log2\bigl(2\log(k+2)+\log(k+3)+3\bigr).
 \label{eq:R-increment-bound}
\end{align}
Define, for real $x\geq2440$,
\begin{align*}
 M(x)=&\ a_1(2x-1)-\frac{\log2}{2}U_{\max}-3\alpha_+\\
 &-S^*\log2\bigl(2\log(x+2)+\log(x+3)+3\bigr).
\end{align*}
Then
\[
 M'(x)=2a_1-S^*\log2\left(\frac2{x+2}+\frac1{x+3}\right)>0,
\]
because the certified lower bound for $d$ gives $2a_1>0.1312$, whereas the
subtracted term is $<0.006$ for $x\geq2440$. Direct substitution gives
$M(2440)>146$. Thus $M(k)>0$ for every $k\geq2440$. Together with
\eqref{eq:R-increment-bound}, this yields
$\varphi(k+1)-\varphi(k)>0$.
\end{proof}

\begin{proof}[Proof of Theorem~\ref{thm:main-full}]
By Lemmas~\ref{lem:phi-base} and \ref{lem:phi-monotone},
$\varphi(k)>0$ for every $k\geq2440$. Proposition~\ref{prop:criterion-phi}
therefore gives $f_k<f_{k+1}$ on that range. The certified ranges in
\eqref{eq:full-sweep-ranges} cover every integer from $6$ through $2439$.
Together the analytic and certified ranges cover all $k\geq6$.
\end{proof}

\section{A general principle}\label{sec:general-principle}

\begin{conjecture}\label{conj:general}
Let $\mathcal P$ be a set of primes with least element $q$ and thin
complement,
\begin{equation}\label{eq:thin-complement}
  \sum_{p\notin\mathcal P}\frac{\log p}{p}<\infty,
\end{equation}
and let $f_k^{\mathcal P}$ be the sum of $1/(n\log n)$ over integers $n$ with
$\Omega(n)=k$ and all prime factors in $\mathcal P$. The sequence
$f_k^{\mathcal P}$ is eventually increasing if and only if $q=2$.
\end{conjecture}

\emph{Some} hypothesis is needed. For $\mathcal P=\{2\}$ one has
$f_k^{\mathcal P}=1/(2^kk\log2)$, and for $\mathcal P=\{2,3\}$ the sequence
again decreases to $0$; in both cases $q=2$ and the conclusion fails. What
goes wrong is that a sparse $\mathcal P$ leaves $\sum_{p\in\mathcal P}p^{-s}$
bounded as $s\to1^+$, so the generating function has no $\zeta(s)^z$ factor
and $f_k^{\mathcal P}$ decays geometrically whatever $q$ is.

Those examples, however, are already excluded by the weaker condition
$\sum_{p\notin\mathcal P}1/p<\infty$, and that condition is what governs the
size of the limit. The candidate limit is
$\prod_{p\notin\mathcal P}(1-1/p)$, which is the value
$H_{\mathcal P}(1,1)$ of \eqref{eq:HP-def} below, and it is a convergent
product with positive value precisely when
$\sum_{p\notin\mathcal P}1/p<\infty$, and $0$ otherwise. That it is indeed the
limit is proved below for $q\geq3$, in
Proposition~\ref{prop:general-interior}(iv); for $\mathcal P$ all primes it is
Lichtman's $f_k\to1$ \cite{Lichtman19}; and for the two sparse examples the
product is $0$, as are their limits. So the unweighted condition separates a
sequence with a positive limit, whose monotonicity is then decided by a
correction term, from one that decays to $0$. The sparse examples are the
degenerate case of a single product, and not a separate phenomenon.

Why then does \eqref{eq:thin-complement} carry the weight $\log p$? Of the two
products in \eqref{eq:HP-def} it is the second, over the complement, that
needs a hypothesis at all. The first converges locally uniformly on
$\sigma>1/2$ for any set of primes whatever, since $\log L_p=O(p^{-2s})$,
while the logarithm of the second is $\asymp-zp^{-s}$ and so requires
$\sum_{p\in\mathcal Q}p^{-\sigma}<\infty$. That is the unweighted condition of
the paragraph above. The weight is needed instead for the $s$-derivative at
$s=1$, where the constant multiplying $2^{-k}$ lives: that constant involves
$\sum_{p\notin\mathcal P}\frac{\log p}{p-1}$, finite precisely
under~\eqref{eq:thin-complement}, and the weaker condition does not give it.
Take the complement to be blocks of primes in $[X_j,X_j^{1+\eta_j}]$ with
$\eta_j\to0$ slowly and $\log X_j$ growing geometrically; then
$\sum_{p\notin\mathcal P}1/p$ converges while
$\sum_{p\notin\mathcal P}(\log p)/p$ diverges. We have not verified that such
a $\mathcal P$ violates the conclusion. We state the hypothesis in the safe
direction, where strengthening it makes the conjecture claim less.

The singularity structure is the one the two theorems above exploit. On
$\abs z<q^\sigma$ the product $\prod_{p\in\mathcal P}(1-zp^{-s})^{-1}$
converges directly, by \eqref{eq:reg-factor-log}, so its first
$z$-singularity is at
$z=q^s$; the factorization
\[
 \prod_{p\in\mathcal P}(1-zp^{-s})^{-1}
  =\prod_p(1-zp^{-s})^{-1}\prod_{p\notin\mathcal P}(1-zp^{-s})
\]
then relates it to $\zeta(s)^zH_{\mathcal P}(s,z)$ on the smaller disc
$\abs z<2^\sigma$ where the middle product also converges. At $s=1$ the point
$z=2$ lies on the singular circle exactly when $q=2$. Note that $z=2$ is
distinguished here not because $2$ is prime, but because it is where the
$\gamma\e^{-2u}$ term of $J$ in Lemma~\ref{lem:zeta-coordinate} evaluates the
regularized product. If $q\geq3$ the extraction at $z=2$ is interior, and
gives a constant multiple of $2^{-k}$ as in Section~\ref{sec:odd}. If $q=2$ it is at the
boundary, and the mechanism producing a factor $k^2$ is the one identified by
Gorodetsky, Lichtman and Wong \cite{GLW} for
$\mathcal P=\{\text{all primes}\}$.

In the interior case $q\geq3$ the regularized product does supply the control
the extraction needs, and the constant can be identified. Write
\begin{equation}\label{eq:local-reg}
 L_p(s,z)=(1-p^{-s})^z(1-zp^{-s})^{-1},
 \qquad
 \log L_p(s,z)=\sum_{m\geq2}\frac{z^m-z}{m}p^{-ms},
\end{equation}
the second identity being \eqref{eq:reg-factor-log}. Define $H_{\mathcal P}$ by the regularized
product, as \eqref{eq:H-full} and \eqref{eq:H-odd} define $H$ and
$H_{\mathrm o}$:
\begin{equation}\label{eq:HP-def}
  H_{\mathcal P}(s,z)
   =\prod_{p\notin\mathcal P}(1-p^{-s})^{z}
    \prod_{p\in\mathcal P}(1-zp^{-s})^{-1}(1-p^{-s})^{z}
   =\prod_{p\notin\mathcal P}(1-p^{-s})^{z}
    \prod_{p\in\mathcal P}L_p(s,z),
\end{equation}
which is \eqref{eq:H-full} when $\mathcal P$ is all primes and
\eqref{eq:H-odd} when $\mathcal P=\{p:p\geq3\}$. For $s>1$, where
$\zeta(s)^{z}$ is finite and non-zero, \eqref{eq:euler-products} extends to
$\mathcal P$ as
\begin{equation}\label{eq:euler-products-P}
  \sum_{k\geq0}P_k^{\mathcal P}(s)z^k
   =\prod_{p\in\mathcal P}(1-zp^{-s})^{-1}
   =\zeta(s)^zH_{\mathcal P}(s,z),
  \qquad
  P_k^{\mathcal P}(s)=\sum_{\substack{\Omega(n)=k\\ p\mid n\Rightarrow p\in\mathcal P}}n^{-s},
\end{equation}
by the same computation, since
$\zeta(s)^{-z}=\prod_p(1-p^{-s})^{z}$ there. It is
\eqref{eq:HP-def} and not the factor $\zeta(s)^{-z}\prod_{\mathcal
P}(1-zp^{-s})^{-1}$ of \eqref{eq:euler-products-P} that is used below: the
latter is $0$ at $s=1$, where every constant in this section is evaluated,
while the former is finite and non-zero there. Put
$h_j^{\mathcal P}(s)=[z^j]H_{\mathcal P}(s,z)$.

\begin{proposition}\label{prop:general-interior}
Let $\mathcal P$ be a set of primes with least element $q\geq3$ and complement
$\mathcal Q$ satisfying \eqref{eq:thin-complement}, and fix
$R\in(2,3)$, which is possible with $R<q$ since $q\geq3$. Then:
\begin{enumerate}
\item[(i)] for every real $s\geq1$ both products in \eqref{eq:HP-def} converge
absolutely and locally uniformly on $\abs z<q^{s}$, where
$H_{\mathcal P}(s,\cdot)$ is therefore analytic and zero-free;
\item[(ii)] $M:=\sup_{s\geq1}\sup_{\abs z=R}\abs{H_{\mathcal P}(s,z)}<\infty$,
and hence $\abs{h_j^{\mathcal P}(s)}\leq MR^{-j}$ for all $j\geq0$ and all
$s\geq1$;
\item[(iii)] each $h_j^{\mathcal P}$ is $C^1$ on $[1,\infty)$ with
$\sup_{s\geq1}\abs{(h_j^{\mathcal P})'(s)}\leq M'R^{-j}$ for some
$M'<\infty$, and
\begin{equation}\label{eq:HP-logderiv}
  \tfrac12\partial_s\log H_{\mathcal P}(1,2)
   =-\sum_{p\in\mathcal P}\frac{\log p}{(p-1)(p-2)}
    +\sum_{p\in\mathcal Q}\frac{\log p}{p-1},
\end{equation}
both sums finite, the second precisely by \eqref{eq:thin-complement};
\item[(iv)] consequently
\begin{equation}\label{eq:general-interior}
  f_k^{\mathcal P}
   =\prod_{p\in\mathcal Q}\Bigl(1-\frac1p\Bigr)
    +c_{\mathcal P}2^{-k}+o(2^{-k}),
  \qquad
  c_{\mathcal P}=H_{\mathcal P}(1,2)
   \bigl(\gamma+\tfrac12\partial_s\log H_{\mathcal P}(1,2)\bigr),
\end{equation}
with
$H_{\mathcal P}(1,2)=\prod_{p\in\mathcal P}(1-p^{-1})^2(1-2p^{-1})^{-1}
 \prod_{p\in\mathcal Q}(1-p^{-1})^2\in(0,\infty)$.
\end{enumerate}
\end{proposition}

\begin{proof}
Parts (i)--(iii) are the analytic groundwork: that $H_{\mathcal P}$ exists,
that it and its $s$-derivative are bounded on a circle $\abs z=R$ with
$R\in(2,3)$, and that the coefficients $h_j^{\mathcal P}$ inherit geometric
decay $R^{-j}$ from that. Part (iv) is the substance, and it uses nothing of
$\mathcal P$ beyond those three facts. The reason the radius must lie strictly
between $2$ and $3$ emerges there: $R>2$ is what makes the coefficient tails
$o(2^{-k})$, and $R<3$ is what makes a geometric sum in the error terms
dominated by its last term.

(i) Fix $s\geq1$ and a compact $K\subset\{\abs z<q^{s}\}$, say
$\abs z\leq R'<q^{s}$. For the product over $\mathcal P$,
\eqref{eq:local-reg} gives $\log L_p(s,z)=O_{R'}(p^{-2s})$ uniformly on $K$ by
\eqref{eq:reg-factor-log}, whose
hypothesis $\abs{zp^{-s}}<1$ holds because $\abs z\leq R'<q^{s}\leq p^{s}$ for
every $p\in\mathcal P$; and $\sum_pp^{-2s}<\infty$ since $s\geq1$. For the
product over $\mathcal Q$, $z\log(1-p^{-s})=O_{R'}(p^{-s})$ and
\[
  \sum_{p\in\mathcal Q}p^{-s}
   \leq\sum_{p\in\mathcal Q}\frac1p
   \leq\frac12+\frac1{\log3}\sum_{p\in\mathcal Q}\frac{\log p}{p}<\infty
  \qquad(s\geq1),
\]
the comparison $1/p\leq(\log p)/(p\log3)$ holding for every $p\geq3$ and the
term $p=2$, which lies in $\mathcal Q$ for every $\mathcal P$ with $q\geq3$
and for which the comparison fails since $\log2<1$, being split off and
bounded by $1/2$. So \eqref{eq:thin-complement} gives
$\sum_{p\in\mathcal Q}1/p<\infty$, which is all that is used here. Both
products therefore converge absolutely and locally uniformly, so each is
analytic in $z$ on $\abs z<q^{s}$; and neither vanishes, because the
convergence of the logarithms exhibits each as the exponential of a finite
quantity.

(ii) The two estimates in (i), applied with $R'=R<3\leq q\leq q^{s}$, depend on
$s$ only through $p^{-s}\leq p^{-1}$, so
\[
  \sup_{s\geq1}\sup_{\abs z=R}\abs{\log H_{\mathcal P}(s,z)}
   \leq C_R\sum_pp^{-2}
    +R\sum_{p\in\mathcal Q}\log\frac1{1-1/p}
   <\infty,
\]
which bounds $M$; the bound on $h_j^{\mathcal P}$ is Cauchy's estimate on
$\abs z=R$.

(iii) Differentiating \eqref{eq:HP-def} logarithmically term by term,
\begin{equation}\label{eq:HP-dlog}
  \partial_s\log H_{\mathcal P}(s,z)
   =-\sum_{p\in\mathcal P}\log p\sum_{m\geq2}(z^m-z)p^{-ms}
    +z\sum_{p\in\mathcal Q}\frac{\log p}{p^{s}-1}.
\end{equation}
For $s\geq1$ and $\abs z\leq R$ the first double sum is dominated termwise by
$C_R(\log p)p^{-2}$ and the second by $R(\log p)/(p-1)$, summable over
$\mathcal Q$ exactly by \eqref{eq:thin-complement}. Both bounds are
independent of $s$, so by the Weierstrass $M$-test the differentiated series
converges uniformly there and the differentiation is legitimate; its sum is
continuous on $[1,\infty)\times\{\abs z\leq R\}$, being a uniformly convergent
series of continuous functions, and
$\partial_sH_{\mathcal P}=H_{\mathcal P}\partial_s\log H_{\mathcal P}$ is
bounded there. Since $[z^j]$ is a Cauchy integral over the circle $\abs z=R$,
which may be held fixed while $s$ varies, differentiation under the integral
sign gives $(h_j^{\mathcal P})'(s)=[z^j]\partial_sH_{\mathcal P}(s,z)$ exactly
as in \eqref{eq:coeff-commute}, and Cauchy's estimate then bounds it by
$M'R^{-j}$. At $z=2$ the inner sum in \eqref{eq:HP-dlog} is
\[
  \sum_{m\geq2}(2^m-2)p^{-m}
   =\frac{4}{p(p-2)}-\frac{2}{p(p-1)}
   =\frac{2}{(p-1)(p-2)},
\]
which gives \eqref{eq:HP-logderiv} after halving.

(iv) The plan is the one of Section~\ref{sec:odd}, carried out with $\mathcal P$ in place
of the odd primes. Pass to the variable $u=\log\zeta(s)$, where the measure
$J(u)\dd u$ has the convergent expansion of
Lemma~\ref{lem:zeta-coordinate}; expand each $h_j^{\mathcal P}$ about $s=1$ in
three groups, its value, its derivative and a remainder; and read off the
first two moments. The value gives the limit, the derivative gives the
$2^{-k}$ coefficient, and everything else has to be shown to be $o(2^{-k})$.

Write $P_k^{\mathcal P}(s)=\sum n^{-s}$ over the $n$ counted by
$f_k^{\mathcal P}$. The proof of
Lemma~\ref{lem:integral-representation} used nothing about the index set beyond
non-negativity of the summands, so it applies verbatim and gives
$f_k^{\mathcal P}=\int_1^\infty P_k^{\mathcal P}(s)\dd s$. Expanding
$\zeta(s)^z=\e^{uz}$ with $u=\log\zeta(s)$ and taking the Cauchy product with
$H_{\mathcal P}(s,z)=\sum_jh_j^{\mathcal P}(s)z^j$ gives
$P_k^{\mathcal P}(s)=\sum_{j=0}^kh_j^{\mathcal P}(s)u^{k-j}/(k-j)!$, which is
\eqref{eq:odd-coefficient-extraction} with $\mathcal P$ in place of the odd
primes; changing variables to $u$ as in \eqref{eq:odd-u-integral}, by
Lemma~\ref{lem:zeta-coordinate}, then gives
\[
  f_k^{\mathcal P}
   =\int_0^\infty\sum_{j=0}^k h_j^{\mathcal P}(s(u))
     \frac{u^{k-j}}{(k-j)!}J(u)\dd u .
\]
Split the integral at $u=1$ and write $s_1=s(1)$, so that $s_1>1$. The point
of splitting there is that the expansions of
Lemma~\ref{lem:zeta-coordinate} are only useful for $u$ bounded away from $0$,
while for small $u$ a crude bound suffices. Indeed the range $u\in(0,1]$ is
$s\geq s_1$, where the contribution is
$\int_{s_1}^\infty P_k^{\mathcal P}(s)\dd s$ and may be bounded without any of the
above, simply because the integers involved are large: every $n$ counted by
$P_k^{\mathcal P}$ has $n\geq q^k\geq3^k$, so
$P_k^{\mathcal P}(s)\leq3^{-k(s-s_1)}P_k^{\mathcal P}(s_1)$ for $s\geq s_1$,
whence
\[
  \int_{s_1}^\infty P_k^{\mathcal P}(s)\dd s
   \leq\frac{P_k^{\mathcal P}(s_1)}{k\log3}
   \leq\frac{2^{-k}}{k\log3}\prod_{p\in\mathcal P}(1-2p^{-s_1})^{-1}
   =O\Bigl(\frac{2^{-k}}{k}\Bigr)=o(2^{-k}),
\]
the middle step because the generating function \eqref{eq:euler-products} for
$\mathcal P$ has non-negative coefficients, so
$P_k^{\mathcal P}(s_1)\leq2^{-k}\sum_jP_j^{\mathcal P}(s_1)2^j$, and the
product converges because $2<3^{s_1}\leq q^{s_1}$. On $[1,\infty)$ the
expansions \eqref{eq:zeta-coordinate} apply with remainders that are
$O(\e^{-3u})$ \emph{uniformly} there: each remainder times $\e^{3u}$ is
continuous on $[1,\infty)$ and bounded as $u\to\infty$ by
Lemma~\ref{lem:zeta-coordinate}, hence bounded on $[1,\infty)$. With
$\int_0^\infty\frac{u^m}{m!}\e^{-au}\dd u=a^{-(m+1)}$,
\begin{equation}\label{eq:gen-moments}
  \int_1^\infty\frac{u^m}{m!}J\dd u=1+\gamma2^{-m}+O(3^{-m}),
  \qquad
  \int_1^\infty\frac{u^m}{m!}\varepsilon J\dd u=\tfrac12 2^{-m}+O(3^{-m}),
\end{equation}
the restriction to $[1,\infty)$ costing at most
$\int_0^1u^m/m!\dd u\leq1/(m+1)!$ in each, which is smaller than the error
terms already present and is absorbed into them throughout.
Now expand the coefficients. Write
\[
  h_j^{\mathcal P}(s(u))=h_j^{\mathcal P}(1)
   +\varepsilon(u)\,(h_j^{\mathcal P})'(1)+r_j(u),
  \qquad
  \abs{r_j(u)}\leq\varepsilon(u)\,\Omega(\varepsilon(u))\,R^{-j},
\]
where
\[
  \Omega(\delta)=\sup_{\abs z=R}\ \sup_{1\leq s\leq1+\delta}
   \abs{\partial_sH_{\mathcal P}(s,z)-\partial_sH_{\mathcal P}(1,z)} .
\]
Here $r_j(u)=\int_1^{s(u)}\bigl((h_j^{\mathcal P})'(t)-(h_j^{\mathcal
P})'(1)\bigr)\dd t$, and Cauchy's estimate on $\abs z=R$ applied to
$\partial_sH_{\mathcal P}(t,\cdot)-\partial_sH_{\mathcal P}(1,\cdot)$ bounds
the integrand by $\Omega(t-1)R^{-j}$, which gives the stated bound on $r_j$
since $\Omega$ is non-decreasing and $s(u)-1=\varepsilon(u)$. Moreover
$\Omega(\delta)\to0$ as $\delta\to0$: by (iii) the function
$\partial_sH_{\mathcal P}$ is continuous on
$[1,\infty)\times\{\abs z\leq R\}$, hence uniformly continuous on the compact
set $[1,2]\times\{\abs z=R\}$, which covers $\delta\leq1$. The three groups now contribute as follows. Two sources of error recur in
each, and both are controlled by the position of $R$: a coefficient tail
$\sum_{j>k}$, which needs $R>2$, and the $O(3^{-m})$ remainders of
\eqref{eq:gen-moments}, which need $R<3$.

\emph{First group, the values $h_j^{\mathcal P}(1)$.} By (ii) the tails
$\sum_{j>k}$ of $\sum_jh_j^{\mathcal P}(1)$ and of
$2^{-k}\sum_jh_j^{\mathcal P}(1)2^{j}$ are $O(R^{-k})$ and
$O(2^{-k}(2/R)^{k})$, both $o(2^{-k})$ since $R>2$. The $O(3^{-m})$ terms
contribute
$O\bigl(3^{-k}\sum_{j\leq k}MR^{-j}3^{j}\bigr)=O(R^{-k})=o(2^{-k})$, where the
middle equality is where $R<3$ is used: it makes the geometric sum dominated
by its last term. So this group gives
$H_{\mathcal P}(1,1)+\gamma H_{\mathcal P}(1,2)2^{-k}+o(2^{-k})$.

\emph{Second group, the derivatives $(h_j^{\mathcal P})'(1)$.} The same
argument with (iii) in place of (ii). The tail of
$2^{-k}\sum_j(h_j^{\mathcal P})'(1)2^{j}$ is $O(2^{-k}(2/R)^{k})$ and the
$O(3^{-m})$ terms are $O(R^{-k})$, so the group gives
$\tfrac12\partial_sH_{\mathcal P}(1,2)2^{-k}+o(2^{-k})$. Here
$\sum_j(h_j^{\mathcal P})'(1)2^{j}$ converges to
$\partial_sH_{\mathcal P}(1,2)$, again because $R>2$.

\emph{Third group, the remainders $r_j$.} This is the only place where
$\Omega$ is needed, and the only $\varepsilon$-argument in the proof. Let
$\delta>0$ and choose $U$ with $\Omega(\varepsilon(U))<\delta$, possible since
$\varepsilon(U)\to0$. On $[U,\infty)$ the terms are at most $\delta R^{-j}$
times the second integral of \eqref{eq:gen-moments}, contributing
$O(\delta2^{-k})$. On $[1,U]$ the integrand is $O_U(1)$ times
$u^{k-j}/(k-j)!$, and $\int_1^Uu^m/m!\dd u\leq U^{m+1}/(m+1)!$, so that range
contributes
$O_U\bigl(\sum_{j\leq k}R^{-j}U^{k-j+1}/(k-j+1)!\bigr)=O_U(R^{-k})=o(2^{-k})$.
As $\delta$ was arbitrary, the group is $o(2^{-k})$.

Collecting the two surviving coefficients gives
$c_{\mathcal P}=\gamma H_{\mathcal P}(1,2)
 +\tfrac12\partial_sH_{\mathcal P}(1,2)$, which is
\eqref{eq:general-interior}. For the limit, $L_p(s,1)=1$, so
$H_{\mathcal P}(s,1)=\prod_{p\in\mathcal Q}(1-p^{-s})$ by \eqref{eq:HP-def}.
The stated product for $H_{\mathcal P}(1,2)$ is \eqref{eq:HP-def} at
$(s,z)=(1,2)$; its factors are positive because $p\geq q\geq3$ on
$\mathcal P$, and it is non-zero by the convergence in (i) at $(s,z)=(1,2)$,
which is legitimate since $2<3\leq q$. Positivity of the factors alone would
not give this.
\end{proof}

By \eqref{eq:general-interior} it is positivity of $c_{\mathcal P}$ that
makes the sequence decrease, and that is now a corollary rather than an open
point.

\begin{corollary}\label{cor:general-positivity}
$c_{\mathcal P}>0$ for every $\mathcal P$ with $q\geq3$
satisfying \eqref{eq:thin-complement}; indeed
$c_{\mathcal P}\geq H_{\mathcal P}(1,2)\,\beta/d$, and
$f_k^{\mathcal P}$ is eventually strictly decreasing.
\end{corollary}

\begin{proof}
Put $a_p=\frac{\log p}{(p-1)(p-2)}\geq0$ and $b_p=\frac{\log p}{p-1}\geq0$ for
$p\geq3$. Since $q\geq3$ we have $2\in\mathcal Q$ and
$\mathcal P=\{p\geq3\}\setminus(\mathcal Q\setminus\{2\})$, with $b_2=\log2$,
so for every admissible $\mathcal Q$, finite or infinite,
\eqref{eq:HP-logderiv} rearranges exactly as
\[
  \tfrac12\partial_s\log H_{\mathcal P}(1,2)
   =\log2-\sum_{p\geq3}a_p
    +\sum_{p\in\mathcal Q\setminus\{2\}}(a_p+b_p),
\]
all three sums converging by \eqref{eq:thin-complement}. The last sum has
non-negative terms, so the right side is least when it is empty, that is at
$\mathcal Q=\{2\}$ and $\mathcal P=\{p:p\geq3\}$, where
$\gamma+\tfrac12\partial_s\log H_{\mathcal P}(1,2)$ equals
$\gamma+\log2-\sum_{p\geq3}\frac{\log p}{(p-1)(p-2)}=\beta/d=0.410970\ldots$,
positive by \eqref{eq:beta-intro}. Since $H_{\mathcal P}(1,2)>0$ by
Proposition~\ref{prop:general-interior}(iv), $c_{\mathcal P}>0$, and
\eqref{eq:general-interior} gives
$f_k^{\mathcal P}-f_{k+1}^{\mathcal P}=c_{\mathcal P}2^{-k-1}+o(2^{-k})>0$ for
all large $k$.
\end{proof}

So the single inequality $\beta>0$ of \eqref{eq:beta-intro} settles the whole
interior case. In particular Corollary~\ref{cor:general-positivity} is one of
the two implications of Conjecture~\ref{conj:general}: it shows that
$f_k^{\mathcal P}$ is eventually strictly decreasing whenever $q\geq3$, so an
admissible $\mathcal P$ with $f_k^{\mathcal P}$ eventually increasing must have
$q=2$. What remains is the converse, and the boundary case $q=2$ lies outside
Proposition~\ref{prop:general-interior}, whose parts (ii), (iii) and (iv) all
require $R>2$ and so an admissible $R<q$; part~(i) needs neither, and
\eqref{eq:HP-def} itself is a definition for any $\mathcal P$ whatever.

The constants of Proposition~\ref{prop:general-interior}(iv) have been checked
numerically, and the two checks should not be conflated. The constant
\emph{formulas} specialize correctly: at $\mathcal Q=\{2\}$ they return
$H_{\mathcal P}(1,2)=0.378695020$ against $d=0.3786950315$ and
$c_{\mathcal P}=0.155632714$ against $\beta$, both to about $10^{-7}$, the
precision of the prime sums involved. The \emph{sequence} $f_k^{\mathcal P}$
approaches its asymptotic far more slowly, and how slowly depends on $q$: at
$\mathcal Q=\{2\}$, where the subleading term is the $O(k^33^{-k})$ of
\eqref{eq:odd-intro-expansion} with constant $4900$, the ratio of
$(f_k^{\mathcal P}-c_0)2^{k}$ to $c_{\mathcal P}$ is still $1.13$ at $k=14$,
the largest index computed, while at $\mathcal Q=\{2,3\}$, where $q=5$, it is
$1.003$. So the $q\geq5$ runs test \eqref{eq:general-interior} end to end and
$\mathcal Q=\{2\}$ calibrates the formulas against a proved case. The
computation is \texttt{genP\_check.py}; alone among the computations here it
uses ordinary floating-point rather than ball arithmetic, carries no
certificate, and nothing else in the paper depends on it.

The boundary case $q=2$ is a different problem, and we do not treat it here.
Separating the factor $L_2$ in \eqref{eq:HP-def} gives the exact identity
$H_{\mathcal P}=H_{\mathcal P\setminus\{2\}}/(1-z2^{-s})$, so the entire
$z$-singularity at $z=2^{s}$ is one explicit pole whose numerator belongs to a
prime set with least element at least $3$. Proposition~\ref{prop:general-interior}
therefore applies to the numerator, and in that sense the crux of the boundary
case is an interior object; the $k^2$ of \eqref{eq:GLW-asymptotic} comes from
the pole. Whether the resulting coefficient depends on that numerator only
through its value at the critical point, which would give
Conjecture~\ref{conj:general}, we have been unable to settle, and we will
return to it separately.

\section{Certification}\label{sec:certification}

The numerical propositions are Proposition~\ref{prop:odd-certified} and
Proposition~\ref{prop:full-certified}. Their statements are finite systems of
interval inequalities, each established in ball arithmetic by Arb/FLINT
\cite{Arb} through its Python bindings, so that every enclosure carries a
rigorous error radius. The scripts, their logs and the certificate are in the
repository cited below, each script recording its settings in its docstring.

What ball arithmetic certifies is the computation. Transcribing an enclosure
into a printed numeral, and combining several into a displayed chain, is
ordinary writing, and that is where the errors in this manuscript have been
found. A late adversarial pass turned up three: a stated bound rounded down
rather than up, so that the text asserted what its own certificate did not
establish; the two-part assembly of $\Sigma_2$, whose displayed chain did not
close under its stated value; and a contour bound attached to the wrong
function, the certified $101.245$ bounding the pole-free part
$G_2(1,z)-\rho/(3-z)$, as every consuming script had always assumed, and not
the regularized product, whose supremum on $\abs z=4.2$ is $228.95$. None
invalidated a conclusion, and each is now checked by script rather than by
eye. \texttt{check\_\allowbreak proof\_\allowbreak spine.py} compares every
printed numeral against the log that certifies it, failing if a stated bound
falls below it; \texttt{verify.py} checks the sweep tables and the $\Gamma$
data against the statements they instantiate, checks the criterion sweeps for
gaps in the union of their ranges, and recomputes the certificate's SHA-256.
One input is not an enclosure but a sum over all primes: $\beta>0$ and the
base case \eqref{eq:odd-base-59} both reduce to an upper bound on
$S=\sum_{p\geq3}(\log p)/((p-1)(p-2))$, certified in
\texttt{odd\_\allowbreak beta\_\allowbreak certified.py} as a partial sum over
$p\leq10^{8}$ plus an integral tail, with $d$ bounded below by a partial
product, every Euler factor of $d$ exceeding $1$ so that no tail is needed.

Two constants are stated above more conservatively than the computation
requires: the threshold $2440$ of Lemma~\ref{lem:phi-base} in place of
$2433$, the first index at which $\varphi$ is positive, and
\eqref{eq:ci-geometric-bound} in place of the sharper bound of
Lemma~\ref{lem:pole-coefficient-bounds}(i). Both sharper values are recorded
in the repository as headroom rather than substituted here, since neither
changes the ranges covered and a printed numeral is only as good as the
certificate behind it.

\begingroup\sloppy
$\Sigma_2$ is not a single computed number but a sum of two, and since it
carries $4395.55$ of the $4857.88$ in the ledger of
Proposition~\ref{prop:odd-remainder} it is worth saying which. The head
$j\leq16$ is a direct interval recursion on the coefficients over $1600$ boxes
in $s$, giving $8423.91$; a Cauchy estimate would be far lossier there, being
dominated by the maximum over the circle and discarding the cancellation among
the true coefficients. The tail $j>16$ is where a Cauchy estimate is cheap:
with $\sup_{\abs z=4.5}$ of the pole-free part at most $1.9672165\cdot10^{6}$,
the tail is at most that times $(3/4.5)^{17}(1-3/4.5)^{-1}$, which is
$5989.94$. The two sum to $14413.85$, so the stated $\Sigma_2\leq14414$ holds
with $0.15$ to spare --- the narrowest margin here except $C_E$. The $40$ boxes
in $s$ behind the tail are coarse enough that it is likely several-fold loose;
that direction is safe, and the budget closes either way.
\endgroup

The constant $C_E$ is certified on the whole of $(1,s(3)]$, and in fact on its
closure. The four quantities entering \eqref{eq:CE-def} are ratios in which
the pole of $\zeta$ at $s=1$ cancels, so each extends continuously to $s=1$.
The Taylor model below is centred at $1+E/2$ with radius $E/2$, and therefore
covers the closed interval $[1,s(3)]$ on which $\Sigma_2$ takes its own
supremum. Interval evaluation on boxes of $s$ loses that cancellation, since
each occurrence of $s-1$ is widened independently, and the enclosures diverge
however finely $(1,s(3)]$ is partitioned. Expanding $\zeta$
as a power series in $\sigma=s-s_0$ about $s_0=1+E/2$, where $E=s(3)-1$, makes
the cancellation symbolic. The resulting series are evaluated on
$\abs{\sigma}\leq E/2$, which is the whole range; their coefficients fall off
geometrically, and truncating after nine terms contributes less than
$10^{-20}$. This is the one constant here with no headroom: the four suprema
are attained at $s(3)$ and give $C_E=5.4890955\ldots$, so the stated
$5.48910$ has about $4\cdot10^{-6}$ to spare, against $0.86\%$ for the ledger.
Any change to the Taylor-model settings must be rechecked against it.

The constant $C_N$ requires $\Pi(t)=\sum_{p\geq3}p^{-t}$ at
$t=s,2s,\dots,200s$ for $s=s(3)$. At $s\approx1.05$ these converge too slowly
to truncate: the Chebyshev-function bounds of Rosser and Schoenfeld
\cite{RS} applied beyond $10^8$ still overestimate the tail twofold,
and that excess compounds through the coefficient recursion. Instead $\Pi$ is
evaluated by the M\"obius series
$\Pi(t)=\sum_{n\geq1}\mu(n)n^{-1}\log\zeta(nt)-2^{-t}$, needing no prime
enumeration; this is the device of Cohen \cite{Cohen}, as used in
\cite[(3.6)]{Lichtman19}. The recursion is run at $600$ bits, cancellation
costing about one digit every $2.3$ steps over the $200$ terms used; at $1200$
bits the head agrees to every printed digit.

For the odd theorem, the coefficient bounds in
\eqref{eq:odd-certified-coefficients} are obtained by cancelling the pole at
$z=3$ before applying coefficient recursion or Cauchy bounds. The range
$0<u\leq3$ is evaluated directly, while the range $u\geq3$ uses the Taylor
model made explicit in Proposition~\ref{prop:odd-two-term}.

For the unrestricted theorem, the $s$-uniform coefficient bounds defining
$\alpha_+$ are obtained by carrying $s$ as an interval over a partition of
$[1,2]$. The signed moments $U$ and $V$ are evaluated before absolute values
are taken, for the reason given at \eqref{eq:UV-exact-combination}.

\subsection*{Availability and versions}
\begingroup\sloppy
Since the contribution is in part a certification, the artefacts and the
versions they were produced with are part of the claim. The scripts, their
logs and the certificate are at
\url{https://github.com/jaideepsaipadhi/erdos-sums-almost-primes}. The
numerics use python-flint~0.9.0; the certificate's \texttt{versions} block
records that binding's version but not the FLINT release under it, which
should be added before submission. The certificate records the parameters of
the \emph{quadrature} run, namely $k_{\max}$, the precision, the panel
geometry and the ellipse, together with a SHA-256 hash of the numerical
payload, which \texttt{verify.py} recomputes. The coefficient sweeps'
parameters are recorded separately, in \texttt{uv\_constants.json} and
\texttt{gamma\_constant.json}, not in the certificate.
\endgroup

\appendix
\section{Certified quadrature}\label{app:quadrature}

Since $1/(n\log n)=\int_1^\infty n^{-s}\dd s$, we have
$f_k=\int_1^\infty P_k(s)\dd s$, and the substitution $s=1+\e^{-u}$ turns this
into
\begin{equation}\label{eq:fk-as-u-integral}
 f_k=\int_{-\infty}^{\infty}F(u)\dd u,
 \qquad F(u)=P_k(1+\e^{-u})\e^{-u},
\end{equation}
the integrand being analytic on the strip $\abs{\operatorname{Im}u}<\pi/2$.
The integral is split into a truncated range $[u_-,u_+]$ with
$u_-=-5.10$ and $u_+=3.5k_{\max}$, taken once for the largest $k$ of the run
and so uniform in $k$, and the two tails beyond it. The truncated range is
divided into panels of width $H=0.8$, each mapped affinely to $[-1,1]$ and
integrated by Gauss--Legendre quadrature with $N=52$ nodes at $1400$ bits of
precision; the number of panels therefore grows linearly with $k_{\max}$.

Four of the six error sources in the certificate are arithmetic: ball
arithmetic throughout, the explicit truncation balls in the M\"obius series
\cite{Cohen} for the prime zeta function, the discretization bound derived below, and
interval-Newton enclosures for the Gauss nodes and weights. The remaining two
are the tails of \eqref{eq:fk-as-u-integral}, which are bounded as follows and
are not negligible by inspection.

For the left tail, $u\leq u_-$ corresponds to $s\geq s_-=1+\e^{5.10}=165.02$.
There $P_1(s)=\sum_pp^{-s}\leq2.78\cdot2^{-s}$, since
$P_1(s)2^{s}=1+\sum_{p\geq3}(2/p)^{s}$ decreases in $s$ and is already below
$2.78$ at $s=2$; and $P_k(s)\leq P_1(s)^k$, because expanding $P_1(s)^k$ counts
every $n$ with $\Omega(n)=k$ at least once. Hence
$P_k(s)\leq2.78^k2^{-ks}$, and
\begin{equation}\label{eq:tail-large-s}
 \int_{-\infty}^{u_-}\!\abs{F}\dd u=\int_{s_-}^\infty\! P_k(s)\dd s
 \leq\frac{2.78^k2^{-ks_-}}{k\log2}.
\end{equation}
For the right tail, $P_k(s)\leq(\log\zeta(s))^k/k!$ for $s>1$
\cite[Lemma 5.6]{Lichtman19}. Substituting $u=\log\zeta(s)$ gives
$\dd s=J(u)\dd u$ with $J=\varepsilon(1+\gamma\varepsilon+\cdots)$ by
Lemma~\ref{lem:zeta-coordinate}, and $\varepsilon=1/\zeta\cdot
(1+\gamma\varepsilon+\cdots)\geq\e^{-u}$, so $J\geq\e^{-u}$ rather than
$\leq$: the clean form $\int_{u_+}^\infty u^k\e^{-u}\dd u/k!$ is not immediate.
What is true is $J(u)\leq\e^{-u}(1+C\e^{-u})$ for $u\geq1$ with $C$ absolute,
and at $u\geq u_+=1750$ the correction is $1+O(\e^{-1750})$, absorbed by the
factor $2$ below with room to spare. With that, the integral is dominated by
$2\int_{u_+}^\infty u^k\e^{-u}\dd u/k!$ and the resulting
incomplete-Gamma integral obeys
$\Gamma(k+1,u_+)/k!\leq2u_+^k\e^{-u_+}/k!$ provided $u_+\geq2k+2$
\cite[\S8.11]{DLMF}, so
\begin{equation}\label{eq:tail-near-1}
 \int_{u_+}^{\infty}\!\abs{F}\dd u\leq\frac{2u_+^k\e^{-u_+}}{k!}.
\end{equation}
The hypothesis $u_+\geq2k+2$ holds throughout, since $u_+=3.5k_{\max}$ and
$k\leq k_{\max}$. Both bounds are added to every reported enclosure. They are
worst at opposite ends of the range, and at their worst they are far below the
differences $f_{k+1}-f_k\asymp a_1k^22^{-k}$ that the sweep has to resolve.
The right-hand side of \eqref{eq:tail-near-1} increases in $k$, since
$u_+$ is fixed and $u_+^k/k!$ is increasing while $k\leq u_+$; that of
\eqref{eq:tail-large-s} is $(2.78\cdot2^{-s_-})^k/(k\log2)$ and so
\emph{decreases}, very fast. The run of record is the one the committed
certificate holds, with $k_{\max}=500$ and so $u_+=1750$ over $2194$ panels.
There \eqref{eq:tail-near-1} is below $10^{-272}$ at its worst, $k=500$, and
\eqref{eq:tail-large-s} below $2.5\cdot10^{-99}$ at its worst, $k=2$, while the
smallest difference in the range actually used, $6\leq k\leq199$, exceeds
$10^{-57}$.

For $f$ analytic inside the Bernstein ellipse $E_\rho$ and bounded there by
$M$, the $N$-node Gauss rule satisfies
\begin{equation}\label{eq:trefethen}
 \left|\int_{-1}^1f(x)\dd x-\widehat I_N\right|
 \leq\frac{64}{15}\,\frac{M\rho^{-2N}}{1-\rho^{-2}}
\end{equation}
\cite[Thm.~4.5]{TrefCC}, where the rule is indexed so that $I_n$ uses $n+1$
nodes and the stated exponent $\rho_B^{-(2n+2)}$ is the $\rho_B^{-2N}$ above;
the bound appears as \cite[Thm.~19.3]{ATAP}. Here $\rho_B$ is the
Bernstein-ellipse parameter, unrelated to the residue $\rho$ of
\eqref{eq:contour-certified}. The implementation evaluates the right side of
\eqref{eq:trefethen} with the interval-scaling factor $H/2$ and a further
factor $10$ of margin. The ellipse is taken with semi-minor axis
$\delta=1.45<\pi/2$ in the $u$ variable, giving $\rho_B\approx7.385$, and $M$ is
bounded using $\abs{P_k(s)}\leq P_k(\operatorname{Re}s)$ together with the
monotonicity of $P_k$ on the reals.

The per-panel bound is not monotone in the panel's position, and it would be
wrong to locate its maximum by inspection: on the right-hand side
$F(u)\approx u^k\e^{-u}/k!$ by \eqref{eq:fk-as-u-integral}, which increases
until $u\approx k$ and only then decays, so $M$ rises across the bulk of the
range before falling. What is used is therefore the maximum over all panels and
all $k\leq k_{\max}$, which the certificate records; the largest radius it
reports is $4.6\times10^{-87}$, attained at $k=2$, and so is negligible beside
the differences the sweep resolves. It also dominates both tails everywhere in
the range: by some $185$ orders of magnitude against
\eqref{eq:tail-near-1}, and by $12$ at $k=2$ against
\eqref{eq:tail-large-s}, where that tail is worst. So the
discretization bound and not the truncation is what limits the accuracy of the
$f_k$.

The pinned coefficient sweeps use primes up to $10^5$ and $300$ bits of
precision, with $1600$ boxes in $s$ for $\alpha_+$ and $800$ boxes for the
second-derivative sweep. These parameters are recorded in
\texttt{uv\_constants.json} and \texttt{gamma\_constant.json} respectively, not
in \texttt{certificate.json}, whose \texttt{parameters} block covers the
quadrature run alone.

\section*{Acknowledgements}

I am grateful to Jared Duker Lichtman for his guidance throughout this work,
and for many helpful discussions.

\section*{Statement on the use of AI}

During the preparation of this work the author used Claude Opus~5 (Anthropic),
a large language model, in order to develop and check proofs, write and run
the interval-arithmetic code, search and check the literature, and draft and
edit the text. The proofs were checked by additional
adversarial review passes of the model, and the numerical claims were
certified in ball arithmetic by scripts available in the repository cited in
Section~\ref{sec:certification}. After using this tool, the author reviewed and edited
the content as needed and takes full responsibility for the content of the
article.

\end{document}